\documentclass{article}
\usepackage{amsmath,amssymb,latexsym}
\usepackage{mathtools}
\usepackage{graphicx}
\usepackage{epsfig}
\usepackage{longtable}
\usepackage{amsthm}
\usepackage{xcolor}
\usepackage{multirow}
\usepackage{caption}
\usepackage{subcaption}
\usepackage{ upgreek }
\usepackage{amssymb}
\usepackage{dsfont}
\usepackage{mathrsfs}
\usepackage{lineno}
\newtheorem{theorem}{Theorem}

\newtheorem{corollary}[theorem]{Corollary}

\newtheorem{definition}{Definition}
\newtheorem{example}{Example}
\newtheorem{remark}{Remark}

\newtheorem{condition}{Condition}

\def\XX{{\mathbf X}}
\def\UU{{\mathbf U}}
\def\YY{{\mathbf Y}}
\def\xx{{\mathbf x}}

\def\uu{{\mathbf u}}
\def\aa{{\mathbf a}}
\def\bb{{\mathbf b}}% Required for inserting images

\title{Directional footrule-coefficients}
\author{Enrique de Amo$^{\rm a}$, David García-Fernández$^{\rm b}$, Manuel Úbeda-Flores$^{\rm a,b,}$\footnote{Corresponding author}\bigskip\\
\small{$^{\rm a}$Department of Mathematics, University of Almería, 04120 Almería, Spain}\\
\small{\texttt{edeamo@ual.es,\,\,\,mubeda@ual.es}}\\
\small{$^{\rm b}$Research Group of Theory of Copulas and Applications, University of Almería, 04120 Almería, Spain}\\
\small{\texttt{dgf992@inlumine.ual.es}}\\
}
\begin{document}

\maketitle

\begin{abstract}
    Rank-based dependence measures such as Spearman’s footrule are robust and invariant, but they often fail to capture directional or asymmetric dependence in multivariate settings. This paper introduces a new family of directional Spearman footrule coefficients for multivariate data, defined within the copula framework to clearly separate marginal behavior from dependence structure. We establish their main theoretical properties, showing full consistency with the classical footrule, including behavior under independence and extreme dependence, as well as symmetry and reflection properties. Nonparametric rank-based estimators are proposed and their asymptotic consistency is discussed. Explicit expressions for several known families of copulas illustrate the ability of the proposed coefficients to detect directional dependence patterns undetected by classical measures.
\end{abstract}
\bigskip

\bigskip


\section{Introduction}

Dependence among random variables has long been a central topic in the statistical literature, owing to its fundamental role in the modeling and analysis of multivariate stochastic systems (see for instance \cite{Jodgeo82}). 

Alternatively, copulas provide a powerful framework for modeling the dependence structure among random variables (\cite{Ahn2015,Durante2016book,Fisher97,Ne06}). By separating the marginal distributions from the joint dependence, they allow for a more nuanced characterization of relationships, capturing features such as tail dependence, asymmetries, and directional dependence that traditional correlation measures cannot. This flexibility makes copulas particularly valuable in multivariate contexts, where understanding the intricate interactions among variables is essential.

In this paper, we investigate the measurement of dependence among the components of a $d$-dimensional random vector ($d\ge 2$) modified along a prescribed direction $\alpha\in\{-1,1\}^d$, within the framework of copula-based methods. More precisely, we introduce and analyse a directional extension of the multivariate Spearman's footrule coefficient. The classical bivariate version of the Spearman's footrule was originally proposed by Spearman in \cite{Sp06}  and later extended to the multivariate setting by \'Ubeda-Flores in \cite{Ub05} (for additional studies of this coefficient, see, e.g., \cite{Chen,Genest,Pe2023}). Our construction is based on multivariate dependence measures derived from the diagonal section of the associated copula $C$. This directional perspective allows for the identification of asymmetric dependence patterns that are not captured by traditional, fully symmetric measures of association.

We show that the extreme cases of the proposed directional coefficients, corresponding to the directions $(1,1,\ldots,1)$ and $(-1,-1,\ldots,-1)$, are consistent with previously introduced multivariate versions available in the literature (see \cite{Dec25}). Moreover, we establish several fundamental properties of these coefficients, highlighting their theoretical coherence and interpretability.

The remainder of the paper is organized as follows. Section~\ref{sec:prel} reviews the necessary background in copula theory and introduces directional dependence concepts that are essential for the development of our results. Section~\ref{sec:main} is devoted to the definition of multidimensional directional footrule-coefficients for arbitrary directions $\alpha$. We derive their main properties and examine their theoretical implications, with particular emphasis on their representation as linear combinations of lower-dimensional directional footrule coefficients. This representation sheds light on the contribution of dependence among subsets of components to the overall dependence structure. We also present illustrative examples for well-known families of copulas. Section~\ref{sec:esti} introduces nonparametric rank-based estimators for the proposed coefficients and presents Monte Carlo simulation results, which illustrate their finite-sample performance and confirm the main theoretical properties. Finally, Section \ref{sec:conc} is devoted to conclusions.

\section{Copulas and directional dependence}\label{sec:prel}

Copulas are a widely used tool in probability and statistics literature. Specifically, a $d$-dimensional copula, $C$, is a function $C:[0,1]^d\longrightarrow [0,1]$ satisfying 
\begin{itemize}
    \item [i)] $C(\uu)=0$ whenever one of the components of $\uu\in[0,1]^d$ equals zero and $C(\uu)=u_k$ whenever every coordinate of $\uu$ equals $1$ but $u_k$.
    \item[ii)] For every $\aa=(a_1,a_2,\ldots,a_d)$ and $\bb=(b_1,b_2,\ldots,b_d)\in [0,1]^d$, such that $a_k\leq b_k$ for all $k=1,\ldots,d$,
        $$V_C([\aa,\bb])=\sum_{{\bf c}}{\rm sgn}({\bf c})C({\bf c})\geq 0,$$
         where $[\aa,\bb]$ denotes the $d$-box $\times_{i=1}^d[a_i,b_i]$ and the sum is taken over the vertices ${\bf c}=(c_1,c_2,\ldots,c_d)$ of $[\aa,\bb],$ i.e., each $c_k$ is equal to either $a_k$ or $b_k$, and ${\rm sgn}({\bf c})=1$ if $c_k=a_k$ for an even number of $k$'s, and ${\rm sgn}({\bf c})=-1$ otherwise.
\end{itemize}

$V_C$ is called the $C$-{\it volume}, and represents the probability mass.

Equivalently, a copula is a multivariate distribution function that is defined as $C(\uu)=\mathbb{P}[U_1\leq u_1,\ldots,U_d\leq u_d],$
where $U_i$ is standard uniform in $[0,1]$. The relevance of copula theory lies in Sklar's theorem (see \cite{Ub2017}):

\begin{theorem}[Sklar]
    Let $\XX=(X_1,X_2,\ldots,X_d)$ be a random $d$-vector with joint distribution function $H$ and one-dimensional marginal distributions $F_1,F_2,\ldots,F_d$. Then there exists a $d$-copula, $C$, uniquely determined on $\times_{i=1}^d Range(F_i)$, such that 
    $$H(\xx)=C(F_1(x_1), F_2(x_2),\ldots,F_d(x_d))$$
    for all $\xx\in[-\infty,+\infty]^d$. If all the marginals $F_i$ are continuous then the $d$-copula $C$ is unique.
\end{theorem}

We refer to \cite{Sk59} for a complete proof.

Independence is represented by the product copula defined as $\Pi_d(\uu)=\prod_{i=1}^du_i$ for all ${\bf u}\in[0,1]^d$, and maximal positive dependence by the upper Fréchet-Hoeffding bound $M_d(\uu)=\min_{1\leq i\leq d}\{u_i\}$, where subscripts in $\Pi_d$ and $M_d$ mean dimension.

Let $\widehat{C}$ be the {\it survival d-copula} of $C$ \cite{Wolff1980}, i.e., $\widehat{C}({\bf u})=\mathbb{P}[{\bf U}>{\bf 1-u}]=\overline{C}(1-\uu)$, where {\bf U} is a random $d$-vector whose components are uniform $[0,1]$ and with associated $d$-copula $C$; and $\overline{C}(\uu)=\mathbb{P[\UU>\uu]}$ is called the {\it survival function}.

The diagonal section of a $d$-copula $C$ is a function $\delta_C : [0, 1] \to [0, 1]$, defined as $\delta_C(u) = C(u, \dots, u)$ for all $u \in [0, 1]$. %$\delta_C$ is the distribution %function of the largest order statistic of a %random vector $\mathbf{U}$---i.e., $U_{\max} = %\max\{U_1, \dots, U_d\}$---, that is,
%$P(U_{\max} \leq u) = P(U_1 \leq u, \dots, U_d %\leq u) = C(u, \dots, u) = \delta_C(u)$
%for all $u \in [0, 1]$, where $C$ is the copula %associated with the random vector $\mathbf{U}$.

For further information about copulas, see \cite{Durante2016book,Ne06}.
 
Two of the most widely recognized concepts used to study potential dependence relationships among random variables are positive (respectively, negative) upper orthant dependence (PUOD) (respectively, NUOD) and  positive (respectively, negative) lower orthant dependence (PLOD) (respectively, NLOD). Intuitively, PUOD and NUOD measure how a set of random variables $X_1,\ldots,X_d$ deviate from independence in terms of their likelihood of attaining large values simultaneously. In particular, under PUOD, random variables tend to attain large values simultaneously more frequently than independent random variables sharing the same marginal distributions, whereas this joint behavior is less pronounced under NUOD (see, e.g., \cite{Joe1997,Mu06,Ne06,Wei2014} for results on positive dependence properties based on copulas). Nonetheless, dependence among random variables may be described through alternative frameworks. For example, for a given subset $J \subset I = \{1,2,\ldots,d\}$, one may encounter dependence structures in which extreme (large or small) realizations of the variables $\{X_j : j \in J\}$ are consistently linked to opposite extreme realizations of the complementary variables $\{X_j : j \in I \setminus J\}$. This concept was investigated in \cite{Que2012}, which motivated the definition presented below.

\begin{definition}
    [PD$(\alpha)$ dependence concept]
    Let $\XX$ be a continuous random vector with joint distribution function $H$, and let $\alpha=(\alpha_1,\alpha_2,\ldots,\alpha_d)\in\mathbb{R}^d$ such that $|\alpha_i|=1$ for all $i=1,2,\ldots,d$. We say that $\XX$ (or $H$) is (orthant) positive (respectively, negative) dependent according to the direction $\alpha$ ---denoted by PD$(\alpha)$ (respectively, ND$(\alpha)$)--- if, for every $\xx\in[-\infty,+\infty]^d$,
\begin{equation}\label{eq:PDalpha}
\mathbb{P}\left[\bigcap_{i=1}^d(\alpha_i X_i>x_i)\right]\geq \prod_{i=1}^d\mathbb{P}[\alpha_i X_i> x_i]
\end{equation}
(respectively, by reversing the sense of the inequality ``$\ge$'' in \eqref{eq:PDalpha}).
\end{definition}

Following from the definition above, the concept PD$({\bf 1})$ (respectively, PD$({\bf -1})$), where ${\bf 1}=(1,1,\ldots,1)$ corresponds to the PUOD (respectively, PLOD); and similarly for the corresponding negative concepts. For more details, see \cite{Que2012,Que2024}.

Based on the PD$(\alpha)$ dependence concept, a dependence (partial) order is provided in \cite{Edeamo2024} with the aim of comparing the strength of the positive dependence in a particular direction of two random vectors.

\begin{definition}[PD$(\alpha)$ order]\label{def:PDorder}
Let $\XX$ and $\YY$ be two random $d$-vectors with respective distribution functions $F$ and $G$. Let $\alpha\in\mathbb{R}^d$ such that $|\alpha_i|=1$ for all $i=1,2,\ldots,d$. $\XX$ is said to be {\it smaller than} $\YY$ {\it in the positive dependence order according to the direction} $\alpha$, denoted by $\XX\leq_{{\rm PD}(\alpha)}\YY$, if, for every $\xx\in\mathbb{R}^n$, we have 
$$\mathbb{P}\left[\bigcap_{i=1}^d(\alpha_iX_i>x_i)\right]\leq\mathbb{P}\left[\bigcap_{i=1}^d(\alpha_iY_i>x_i)\right].$$
\end{definition}

Let $\XX$ be a $d$-dimensional random vector with associated $d$-copula $C$ and marginals uniformly distributed on $[0,1]$. Let ${\cal{J}}=\{i_1,i_2,\ldots,i_m\}\subseteq\{1,2,\ldots,d\}$ and ${\cal{I}}=\{1,2,\ldots,d\}\setminus{\cal{J}}$ such that $i_1<i_2<\cdots<i_m$. Let $C_{\cal{J}}:[0,1]^m\longrightarrow [0,1]$ be the ${\cal{J}}$-{\it marginal} copula of $C$ for the random vector $(X_{i_1},X_{i_2},\ldots,X_{i_m})$, which is defined by setting $n-m$ arguments of $C$ equal to $1$, i.e.,
$$C_{\cal{J}}(u_1,u_2,\ldots,u_m)=C(v_1,v_2,\ldots,v_n),$$
where $v_i=u_i$ if $i\in{\cal{J}}$, and $v_i=1$ otherwise. For any subset ${\cal{J}} \subseteq \{1,2,...,d\}$, we define the complemented random vector
\begin{equation}\label{eq:Xm}
{\bf X}^{({\cal{J}})} = (Y_1,\dots,Y_n), 
\qquad
Y_j =
\begin{cases}
X_j, & j \notin {\cal{J}} ,\\[4pt]
1-X_j, & j \in {\cal{J}} ,
\end{cases}
\end{equation}
and denote by $C^{({\cal{J}} )}$ its associated $d$--copula. The number of complemented coordinates is $m = |{\cal{J}} |$. Its associated copula can be expressed as follows:

\begin{theorem}\label{th:main}
For any ${\cal{J}} \subseteq \{1,2,\ldots,d\}$, the copula $C^{({\cal{J}})}$ associated with ${\bf X}^{({\cal{J}})}$ is given by
\begin{equation}\label{eq:main}
C^{({\cal{J}})}({\bf u})=\sum_{\mathcal{K}\subseteq \mathcal{J}}(-1)^{|\mathcal{K}|}\, C_{\mathcal{I}\cup \mathcal{K}}\Bigl(\nu({\bf u},\, \mathcal{I}, \mathcal{K})\Bigr)
\end{equation}
for all ${\bf u}\in[0,1]^d$, where, for each subset $\mathcal{K}\subseteq \mathcal{J}$, the mapping $\nu({\bf u},\, \mathcal{I}, \mathcal{K})$ is defined such that it assigns $u_i$ if $i\in \mathcal{I}$, and $1-u_i$ if $i\in \mathcal{K}$.
\end{theorem}

\begin{proof}Let $\mathbf{X}^{({\cal{J}})}=(Y_1,\ldots,Y_d)$ be the random vector whose components are given by \eqref{eq:Xm}. For $j\notin\mathcal{J}$ the inequality $Y_j\le u_j$ is equivalent to $X_j\le u_j$, and for $j\in\mathcal{J}$ we have $X_j\ge 1-u_j$. Thus, the $d$-copula $C^{({\cal{J}})}$ associated with $\mathbf{X}^{({\cal{J}})}$ is given by
$$C^{({\cal{J}})}({\bf u})=\mathbb{P}\left[\mathbf{X}^{({\cal{J}})}\le{\bf u}\right]=\mathbb{P}\left[\bigcap_{j\notin\mathcal{J}}\{X_j\le u_j\},\bigcap_{j\in\mathcal{J}}\{X_j\ge 1-u_j\}\right]$$
for all ${\bf u}\in[0,1]^d$. Define
$$A=\bigcap_{j\notin\mathcal{J}}\{X_j\le u_j\}.$$
Since for $j\in\mathcal{J}$ the event $\{X_j\ge 1-u_j\}$ is the complement of $\{X_j\le 1-u_j\}$, we have
$$\begin{aligned}
\left(\bigcap_{j\notin\mathcal{J}}\{X_j\le u_j\},\bigcap_{j\in\mathcal{J}}\{X_j\ge 1-u_j\}\right)& = \Biggl\{ \omega\in\Omega:\, \omega\in A\text{ and } \omega\notin \bigcup_{j\in\mathcal{J}}\{X_j\le 1-u_j\}\Biggr\}\\
& = A\setminus \left(A\cap\bigcup_{j\in\mathcal{J}}\{X_j\le 1-u_j\}\right).
\end{aligned}$$
Thus,
\begin{equation}\label{eq:expC}
C^{({\cal{J}})}({\bf u})=\mathbb{P}[A]-\mathbb{P}\left[A\cap\bigcup_{j\in\mathcal{J}}\{X_j\le 1-u_j\}\right].
\end{equation}
By the inclusion--exclusion principle, we can express the probability of the union as
$$\mathbb{P}\left[A\cap\bigcup_{j\in\mathcal{J}}\{X_j\le 1-u_j\}\right]
=\sum_{\emptyset\neq \mathcal{K}\subseteq \mathcal{J}}(-1)^{|\mathcal{K}|-1}\,\mathbb{P}\left[A\cap\bigcap_{j\in K}\{X_j\le 1-u_j\}\right].$$
Since $C$ is the copula of $\mathbf{X}$, for any subset $\mathcal{L}\subseteq \{1,2,\ldots,d\}$ and any numbers $\{t_j\}_{j\in \mathcal{L}}\subseteq[0,1]$ we have
$$\mathbb{P}\left[\bigcap_{j\in \mathcal{L}}\{X_j\le t_j\}\right]=C_\mathcal{L}\Bigl((t_j)_{j\in \mathcal{L}}\Bigr).$$
Then we have
$$\mathbb{P}(A)=\mathbb{P}\left[\bigcap_{j\notin\mathcal{J}}\{X_j\le u_j\}\right]
=C_{\mathcal{I}}\Bigl(\nu(\mathbf{u},\mathcal{I})\Bigr),$$
and for any nonempty $\mathcal{K}\subseteq\mathcal{J}$,
$$\mathbb{P}\left[A\cap\bigcap_{j\in \mathcal{K}}\{X_j\le 1-u_j\}\right]
=C_{\mathcal{I}\cup \mathcal{K}}\Bigl(\nu(\mathbf{u},\mathcal{I}, \mathcal{K})\Bigr).$$
Substituting the above expressions Eq. \eqref{eq:expC}, we obtain
$$\begin{aligned}
C^{({\cal{J}})}({\bf u}) & = C_{\mathcal{I}}\Bigl(\nu(\mathbf{u},\mathcal{I})\Bigr)
-\sum_{\emptyset\neq \mathcal{K}\subseteq \mathcal{J}}(-1)^{|\mathcal{K}|-1}\, C_{\mathcal{I}\cup \mathcal{K}}\Bigl(\nu(\mathbf{u},\mathcal{I}, \mathcal{K})\Bigr)\\[1mm]
& = C_{\mathcal{I}}\Bigl(\nu(\mathbf{u},\mathcal{I})\Bigr)
-\sum_{1\le k\le m} C_{\mathcal{I}\cup\{i_k\}}\Bigl(\nu(\mathbf{u},\mathcal{I},\{i_k\})\Bigr)\\[1mm]
&\quad +\sum_{1\le k<l\le m} C_{\mathcal{I}\cup\{i_k,i_l\}}\Bigl(\nu(\mathbf{u},\mathcal{I},\{i_k,i_l\})\Bigr)+\cdots+(-1)^m\,C\Bigl(\nu(\mathbf{u},\mathcal{I},\mathcal{J})\Bigr),
\end{aligned}$$
where we have enumerated the elements of $\mathcal{J}$ as $\{i_1,i_2,\ldots,i_m\}$. Therefore, Expression \eqref{eq:main} easily follows.
\end{proof}

\begin{remark}We want to stress that, while Theorem \ref{th:main} formally appears in \cite[Theorem 4.1]{Fuc14} with a proof conducted via induction, we provide an alternative proof. We believe this different methodological approach facilitates a clearer conceptual understanding of our overall framework in this paper.
\end{remark}
%\begin{remark}
%    The independence copula $\Pi$ is invariant %under any decreasing transformation of its %components.
%\end{remark}

\section{Directional version of multivariate Spearman's footrule}\label{sec:main}

Numerous coefficients have been proposed to capture positive or negative dependence among random variables (see, e.g., \cite{Dol06,Joe1990,Ne96}); however, those designed to measure the directional dependence of the components of a random vector are considerably less common, particularly in the multivariate context---for dimensions 2 and 3, see, e.g., \cite{Gar2013}. Nevertheless, they generalize the former and are capable of capturing dependence relationships that the traditional coefficients fail to detect.
Inspired by the multivariate version of Spearman's footrule for a $d$-dimensional random vector ${\bf U}$ with copula $C$ introduced in \cite{Ub05} as 
\begin{equation*}\label{eq:Ub}\varphi_d(C)=\frac{d+1}{d-1}\int_0^1\left[\delta_C(u)-\delta_{\widehat{C}}(u)\right]{\rm d}u-\frac{2}{d-1},
\end{equation*}
in \cite{Dec20} two multivariate versions of this coefficient were described, called {\it downward} and {\it upward diagonal dependence indices}, respectively:
$$\varphi_d^-(C)=\frac{\int_0^1\left[\delta_C(u)-\delta_\Pi(u)\right] {\rm d}u}{\int_0^1\left[\delta_M(u)-\delta_\Pi(u) \right]{\rm d}u}=\frac{2(d+1)}{(d-1)}\int_0^1\delta(u){\rm d}u-\frac{2}{d-1}$$
and
$$\varphi_d^+(C)=\frac{\int_0^1\left[\delta_{\widehat{C}}(u)-\delta_{\widehat{\Pi}}(u)\right]{\rm d}u}{\int_0^1\left[\delta_{\widehat{M}}(u)-\delta_{\widehat{\Pi}}(u)\right]{\rm d}u}=\frac{2(d+1)}{(d-1)}\int_0^1\delta_{\widehat{C}}(u){\rm d}u-\frac{2}{d-1}$$
(see also \cite{Dec25}). Intuitively, the first coefficient measures the area between the diagonal section of $C$ and the diagonal section of the independence copula, normalized by the full area between the main diagonal, which represents maximal positive dependence and the curve representing independence, and the second one measures the normalized area between the survival function of the minimal position and the curve that represents independence. It is possible to rewrite the coefficients as 
$$\varphi_d^-(C)=\frac{2(d+1)}{(d-1)}\int_0^1\mathbb{P}[U_1\leq u,\ldots,U_d\leq u]{\rm d}u-\frac{2}{d-1} \quad\text{and}\quad\varphi_d^+(C)=\frac{2(d+1)}{(d-1)}\int_0^1\mathbb{P}[U_1> u,\ldots,U_d> u]{\rm d}u-\frac{2}{d-1}.$$

From these last expressions, we define a {\it directional version of Spearman's footrule} (or {\it directional footrule-coefficients}, for short) as follows
\begin{equation}\label{eq:main2}
    \varphi_d^\alpha(C)=\frac{2(d+1)}{d-1}\int_0^1\left(\mathbb{P}[\alpha_1U_1>\alpha_1u,\ldots,\alpha_dU_d>\alpha_du]-\prod_{i=1}^d\mathbb{P}[\alpha_iU_i>\alpha_iu]\right){\rm d}u.
\end{equation}
When $\alpha={\bf1}$, we have $\varphi_d^\alpha=\varphi^+_d$ and when $\alpha={\bf-1}$, then $\varphi^\alpha_d=\varphi_d^-$. Since the variables are uniformly distributed on $[0,1]$ and  the event $U_i > u$ has probability $1-u$, it holds
$$\varphi_d^\alpha(C)=\frac{2(d+1)}{d-1}\int_0^1\mathbb{P}[\alpha_1U_1> \alpha_1u,...,\alpha_dU_d> \alpha_du]{\rm d}u-\frac{2}{(d-1)\binom{d}{|\mathcal{J}|}}.$$

We can rewrite the probabilities in (\ref{eq:main2}) using indicator functions and expectations. Let $\mathbf{1}(A)$ denote the indicator function of event $A$. Applying this, we get:
$$\varphi_d^\alpha(C)=\frac{2(d+1)}{d-1}\int_0^1\left(\mathbb{E}\left[\prod_{i=1}^d \mathbf{1}(\alpha_iU_i > \alpha_iu)\right]-\prod_{i=1}^d \mathbb{E}[\mathbf{1}(\alpha_iU_i > \alpha_iu)]\right){\rm d}u.$$
Here, the expectation $\mathbb{E}$ is taken with respect to the joint distribution of the random vector ${\bf U}$ under the $d$-copula $C$. Note that, since the integrand is bounded by 1, by Fubini-Tonelli we have $$\int_0^1\mathbb{E}\left[\prod_{i=1}^d\mathbf{1}(\alpha_iU_i>\alpha_iu)\right]{\rm d}u=\mathbb{E}\left[\int_0^1\prod_{i=1}^d\mathbf{1}(\alpha_iU_i>\alpha_iu){\rm d}u\right],$$ and $\prod_{i=1}^d\mathbf{1}(\alpha_iU_i>\alpha_iu)=\mathbf{1}(u<\min_{i\in J}U_i,u>\max_{i\in I}U_i)$; therefore, $$\int_0^1\prod_{i=1}^d\mathbf{1}(\alpha_iU_i>\alpha_iu){\rm d}u=\max\left(\min_{i\in J}U_i-\max_{i\in I}U_i,0\right),$$  
since $\prod_{i=1}^d\mathbf{1}(\alpha_iU_i>\alpha_iu)=1$ if, and only if, $\max_{i\in I}U_i<u<\min_{i\in J}U_i.$ On the other hand, 
$$\int_0^1\prod_{i=1}^d\mathbb{E}\left[\mathbf{1}(\alpha_iU_i>\alpha_iu)\right]du=\frac{|I|!|J|!}{(d+1)!}.$$
Thus, the formula can be written as:
\begin{equation}\label{eq:esp}
    \varphi_d^\alpha(C)=\frac{2(d+1)}{d-1}\left( \mathbb{E}\left[\left(\min_{i\in J}U_i-\max_{i\in I}U_i\right)_+\right]-\frac{|I|!|J|!}{(d+1)!}\right),
\end{equation}
where $(\cdot)_+$ denotes $\max(\cdot,0)$.

Some of the main properties of these directional footrule-coefficients are listed in the following result.

\begin{corollary}
    Let $\UU$ be a $d$-dimensional random vector with uniformly distributed marginals on $[0,1]$ and associated copula $C$, then
    \begin{itemize}
        \item [i)] $\sum_\alpha\varphi_d^\alpha(C)=0$.
        \item [ii)] $\varphi_d=\sum_{\alpha\neq\mathbf{1},-\mathbf{1}}-\frac{1}{2}\varphi_d^\alpha$. %where $\varphi_d$ is the multidimensional Spearman's footrule given by \eqref{eq:Ub}.
        \item[iii)]  $\varphi_d^\alpha(C)=\varphi_d^{-\alpha}(\widehat{C})$.
    \end{itemize}
\end{corollary}

\begin{proof}
    Part $i)$ is straightforward since we have $$\sum_{\alpha}\mathbb{P}\left[\bigcap_{i=1}^d(\alpha_i U_i>\alpha_i u)\right]=\sum_{\alpha}\prod_{i=1}^d\mathbb{P}[\alpha_i U_i>\alpha_iu]=1.$$

    In order to prove $ii)$, starting from $i)$ $\sum_\alpha\varphi_d^\alpha=0$, so we have $\sum_{\alpha\neq\mathbf{1},-\mathbf{1}}\varphi_d^\alpha+\varphi_d^++\varphi_d^-=0$. Then
    $$\frac{\varphi^+_d+\varphi_d^-}{2}=\sum_{\alpha\neq\mathbf{1},-\mathbf{1}}-\frac{1}{2}\varphi_d^\alpha,.$$
    In \cite{Dec25} it is shown that $\varphi_d=(\varphi_d^++\varphi_d^-)/2$ and then, the desired expression is obtained.
    
    For part $iii)$, we know that $\widehat{C}(\uu)=\mathbb{P}[U>1-u]=\mathbb{P}[-U<u-1]=\mathbb{P}[-U>-u]$ since $U$ is uniform in $[0,1]$. Therefore,
    \begin{align*}
        \varphi_d^\alpha(C)&=\frac{2(d+1)}{d-1}\int_0^1\left(\mathbb{P}\left[\bigcap_{i=1}^d(\alpha_iU_i>\alpha_iu)\right]-\prod_{i=1}^d\mathbb{P}[\alpha_iU_i>\alpha_iu]\right){\rm d}u\\
        &=\frac{2(d+1)}{d-1}\int_0^1\left(\mathbb{P}\left[\bigcap_{i=1}^d(-\alpha_i(-U_i)>-\alpha_i(-u))\right]-\prod_{i=1}^d\mathbb{P}[-\alpha_i(-U_i)>-\alpha_i(-u)]\right){\rm d}u\\
        &=\varphi_d^{-\alpha}(\widehat{C}),
    \end{align*}
    as desired.
\end{proof}

In \cite{Ne2011} it was conjectured that multidimensional directional $\rho$-coefficients---another type of multidimensional measure of directionality---can be expressed as a linear combination of lower dimensional directional coefficients. That conjecture was proved in \cite{Edeamo25}. Now we show that the new directional coefficients can also be expressed in terms of lower dimensional coefficients. 

\begin{theorem}\label{th:exp}
    Let $\XX$ be a $d$-dimensional random vector whose marginal are uniform on $[0,1]$, let $C$ be its associated $d$-copula and $\alpha\in\mathbb{R}^d$ such that $\alpha_i\in\{-1,1\}$ for all $i=1,\ldots,d$. Let $I\subseteq\{1,\ldots,d\}$ such that $\alpha_i=-1$ if $i\in I$ and $\alpha_i=1$ if $i\in J=\{1,\ldots,d\}\setminus I$. Then
    \begin{equation}\label{th:expr}
        \varphi_d^\alpha(C)=\frac{2(d+1)}{d-1}\sum_{S\subseteq J}(-1)^{|S|}\frac{|I|+|S|-1}{2(|I|+|S|+1)}\varphi^-_{{\bf X}_{I\cup S}},
    \end{equation}
where ${\bf X}_K$ is the random vector whose components are $X_i$ for $i\in K$, and $\varphi_{{\bf X}_{K}}^-$ corresponds to the coefficient $\varphi_k^-$ for ${\bf X}_K$.
\end{theorem}

\begin{proof}
    By using (\cite[Lemma 3.4 and Remark 3.5]{Edeamo}), it holds
    \begin{align*}
        \varphi_d^\alpha(C)&=\frac{2(d+1)}{d-1}\int_0^1\mathbb{P}[\alpha_1U_1> \alpha_1u,\ldots,\alpha_dU_d> \alpha_du]{\rm d}u-\frac{2}{(d-1)\binom{d}{|\mathcal{J}|}}\\
        &=\frac{2(d+1)}{d-1}\int_0^1\left( \sum_{k=0}^{|J|} (-1)^k \sum_{\substack{S \subseteq J \\ |S|=k}} C(\mathbf{u}_{I \cup S}, \mathbf{1})\right){\rm d}u-\frac{2}{(d-1)\binom{d}{|\mathcal{J}|}}.\\
    \end{align*}
After some elementary calculations, we obtain
$$\int_0^1 C(\mathbf{u}_{I \cup S},{\bf 1}) {\rm d}u=\frac{|I|+|S|-1}{2(|I|+|S|+1)}\left(\varphi_{\XX_{I\cup S}}^-+\frac{2}{|I|+|S|-1}\right)=\frac{|I|+|S|-1}{2(|I|+|S|+1)}\varphi_{\XX_{I\cup S}}^-+\frac{1}{|I|+|S|+1},$$
since $$\varphi_{\XX_{I\cup S}}^-(C)=\frac{2(|I|+|S|+1)}{|I|+|S|-1}\int_0^1C(\uu_{I\cup S},{\bf 1})du-\frac{2}{|I|+|S|-1}.$$ So, using this equality we have
\begin{align*}
    \varphi_d^\alpha(C)&=\frac{2(d+1)}{d-1}\left(\sum_{k=0}^{|J|}(-1)^k\sum_{\substack{S\subseteq J\\ |S|=k}}\left\{\frac{|I|+|S|-1}{2(|I|+|S|+1)}\varphi_{\XX_{I\cup S}}^-+\frac{1}{|I|+|S|+1}\right\}\right)-\frac{2}{(d-1)\binom{d}{|J|}}\\
    &=\frac{2(d+1)}{d-1}\sum_{S\subseteq J}(-1)^{|S|}\frac{|I|+|S|-1}{2(|I|+|S|+1)}\varphi_{\XX_{I\cup S}}^-.
\end{align*}
The terms that do not involve the $\varphi^-$ coefficients vanish because 
$$\sum_{S\subseteq J}(-1)^{|S|}\frac{1}{|I|+|S|+1}=\sum_{k=0}^{|J|}(-1)^k\binom{|J|}{k}\frac{1}{|I|+k+1}=\frac{|I|!|J|!}{(d+1)!}$$
(see, e.g., \cite{Graham}), and $$\frac{2(d+1)}{d-1}\cdot\frac{|I|!|J|!}{(d+1)!}=\frac{2}{(d-1)\binom{d}{|J|}}.$$
Therefore, the desired formula is obtained.
\end{proof}

Table \ref{tab:f_coefficients} shows the sixteen directional footrule-coefficients for $d = 4$  in terms of lower dimensional footrule-
coefficients.
\begin{table}[htb]
\begin{center}
\begin{longtable}{c p{0.8\textwidth}}
\hline
\textbf{Direction $\alpha$} & \textbf{Coefficient $\varphi_4^\alpha(C)$} \\
\hline
\endhead % End of table header for longtable
\hline
\multicolumn{2}{|r|}{{Continued on next page}} \\
\endfoot % Footer for continuing pages
\hline
\endlastfoot % Footer for the last page
$(1,1,1,1)$ & $\frac{5}{9}(\varphi^-_{\mathbf{X}_{12}} + \varphi^-_{\mathbf{X}_{13}} + \varphi^-_{\mathbf{X}_{14}} + \varphi^-_{\mathbf{X}_{23}} + \varphi^-_{\mathbf{X}_{24}} + \varphi^-_{\mathbf{X}_{34}}) - \frac{5}{6}(\varphi^-_{\mathbf{X}_{123}} + \varphi^-_{\mathbf{X}_{124}} + \varphi^-_{\mathbf{X}_{134}} + \varphi^-_{\mathbf{X}_{234}}) + \varphi^-_{\mathbf{X}_{1234}}$ \\

$(-1,1,1,1)$ & $-\frac{5}{9}(\varphi^-_{\mathbf{X}_{12}} + \varphi^-_{\mathbf{X}_{13}} + \varphi^-_{\mathbf{X}_{14}}) + \frac{5}{6}(\varphi^-_{\mathbf{X}_{123}} + \varphi^-_{\mathbf{X}_{124}} + \varphi^-_{\mathbf{X}_{134}}) - \varphi^-_{\mathbf{X}_{1234}}$ \\

$(1,-1,1,1)$ & $-\frac{5}{9}(\varphi^-_{\mathbf{X}_{12}} + \varphi^-_{\mathbf{X}_{23}} + \varphi^-_{\mathbf{X}_{24}}) + \frac{5}{6}(\varphi^-_{\mathbf{X}_{123}} + \varphi^-_{\mathbf{X}_{124}} + \varphi^-_{\mathbf{X}_{234}}) - \varphi^-_{\mathbf{X}_{1234}}$ \\

$(1,1,-1,1)$ & $-\frac{5}{9}(\varphi^-_{\mathbf{X}_{13}} + \varphi^-_{\mathbf{X}_{23}} + \varphi^-_{\mathbf{X}_{34}}) + \frac{5}{6}(\varphi^-_{\mathbf{X}_{123}} + \varphi^-_{\mathbf{X}_{134}} + \varphi^-_{\mathbf{X}_{234}}) - \varphi^-_{\mathbf{X}_{1234}}$ \\

$(1,1,1,-1)$ & $-\frac{5}{9}(\varphi^-_{\mathbf{X}_{14}} + \varphi^-_{\mathbf{X}_{24}} + \varphi^-_{\mathbf{X}_{34}}) + \frac{5}{6}(\varphi^-_{\mathbf{X}_{124}} + \varphi^-_{\mathbf{X}_{134}} + \varphi^-_{\mathbf{X}_{234}}) - \varphi^-_{\mathbf{X}_{1234}}$ \\

$(-1,-1,1,1)$ & $\frac{5}{9}\varphi^-_{\mathbf{X}_{12}} - \frac{5}{6}(\varphi^-_{\mathbf{X}_{123}} + \varphi^-_{\mathbf{X}_{124}}) + \varphi^-_{\mathbf{X}_{1234}}$ \\

$(-1,1,-1,1)$ & $\frac{5}{9}\varphi^-_{\mathbf{X}_{13}} - \frac{5}{6}(\varphi^-_{\mathbf{X}_{123}} + \varphi^-_{\mathbf{X}_{134}}) + \varphi^-_{\mathbf{X}_{1234}}$ \\

$(-1,1,1,-1)$ & $\frac{5}{9}\varphi^-_{\mathbf{X}_{14}} - \frac{5}{6}(\varphi^-_{\mathbf{X}_{124}} + \varphi^-_{\mathbf{X}_{134}}) + \varphi^-_{\mathbf{X}_{1234}}$ \\

$(1,-1,-1,1)$ & $\frac{5}{9}\varphi^-_{\mathbf{X}_{23}} - \frac{5}{6}(\varphi^-_{\mathbf{X}_{123}} + \varphi^-_{\mathbf{X}_{234}}) + \varphi^-_{\mathbf{X}_{1234}}$ \\

$(1,-1,1,-1)$ & $\frac{5}{9}\varphi^-_{\mathbf{X}_{24}} - \frac{5}{6}(\varphi^-_{\mathbf{X}_{124}} + \varphi^-_{\mathbf{X}_{234}}) + \varphi^-_{\mathbf{X}_{1234}}$ \\

$(1,1,-1,-1)$ & $\frac{5}{9}\varphi^-_{\mathbf{X}_{34}} - \frac{5}{6}(\varphi^-_{\mathbf{X}_{134}} + \varphi^-_{\mathbf{X}_{234}}) + \varphi^-_{\mathbf{X}_{1234}}$ \\

$(-1,-1,-1,1)$ & $\frac{5}{6}\varphi^-_{\mathbf{X}_{123}} - \varphi^-_{\mathbf{X}_{1234}}$ \\

$(-1,-1,1,-1)$ & $\frac{5}{6}\varphi^-_{\mathbf{X}_{124}} - \varphi^-_{\mathbf{X}_{1234}}$ \\

$(-1,1,-1,-1)$ & $\frac{5}{6}\varphi^-_{\mathbf{X}_{134}} - \varphi^-_{\mathbf{X}_{1234}}$ \\

$(1,-1,-1,-1)$ & $\frac{5}{6}\varphi^-_{\mathbf{X}_{234}} - \varphi^-_{\mathbf{X}_{1234}}$ \\

$(-1,-1,-1,-1)$ & $\varphi^-_{\mathbf{X}_{1234}}$ \\
\hline
\end{longtable}
\caption{The sixteen directional footrule-coefficients for $d=4$ expressed in terms of lower dimensional footrule-coefficients.} \label{tab:f_coefficients}
\end{center}
\end{table}

Next we provide several examples of directional footrule-coefficients for well-known families of copulas.

\begin{example}[Directional footrule-coefficient for $\Pi_d$]
    If $C=\Pi_n$ then we have $\varphi_d^\alpha(\Pi_d)=0$ since, in this case, $$\mathbb{P}\left[\bigcap_{i=1}^d(\alpha_i X_i>x_i)\right]=\prod_{i=1}^d\mathbb{P}[\alpha_i X_i> x_i].$$
\end{example}

\begin{example}[Directional footrule-coefficients for $M_d$]
For the minimum copula $C = M_d$, the random variables satisfy
$X_1 = X_2 = \cdots = X_d$ in distribution. Since the marginals of the
copula are uniformly distributed on $[0,1]$, each associated random
variable $X_i$ is also uniformly distributed on $[0,1]$. Consequently,
we consider
\begin{equation*}\label{eq:Probab}
v_i(u, \alpha_i) := \mathbb{P}[\alpha_i X_i > \alpha_i u] = \begin{cases} 1 - u & \text{if } \alpha_i = 1, \\ u & \text{if } \alpha_i = -1. \end{cases}
\end{equation*}
We consider the integral
$$I_d(\alpha) = \int_0^1 \left( \mathbb{P}\left[\bigcap_{i=1}^d(\alpha_iX>\alpha_iu)\right] - \prod_{i=1}^d v_i(u, \alpha_i) \right) {\rm d}u$$
for the cases where all $\alpha_i$ are the same. We study several cases depending on $k$, i.e. the number of components of $\alpha$ equal to $-1$.
\begin{enumerate}
\item $\alpha = {\bf 1}$ ($k=0$).
The joint probability is $\mathbb{P}\left[\bigcap_{i=1}^d(X > u)\right] = 1 - u$.
The product of marginal probabilities is $\prod_{i=1}^d (1 - u)$.
The integral becomes:
$$
I_d({\bf 1}) = \int_0^1 \left( 1 - u - \prod_{i=1}^d (1 - u) \right) {\rm d}u= \frac{d-1}{2(d+1)}.
$$
Substituting into the definition of $\varphi_d^\alpha(M_d)$:
$$
\varphi_d^{{\bf 1}}(M_d) = \frac{2(d+1)}{d-1}\cdot\frac{d-1}{2(d+1)} = 1.
$$
\item $\alpha ={\bf -1}$ ($k=d$).
The joint probability is $\mathbb{P}\left[\bigcap_{i=1}^d(X < u)\right] = u$.
The product of marginal probabilities is $\prod_{i=1}^d u$.
The integral becomes:
$$
I_d({\bf -1}) = \int_0^1 \left( u - \prod_{i=1}^d u \right) {\rm d}u= \frac{d-1}{2(d+1)}.
$$
Substituting into the definition of $\varphi_d^\alpha(M_d)$:
$$
\varphi_d^{{\bf -1}}(M_d) = \frac{2(d+1)}{d-1}\cdot\frac{d-1}{2(d+1)} = 1.
$$
\item $k$ is even ($0 < k < d$). When $k$ is even, we have $k$ variables where we are looking at the probability that $-X_i > -u$ (i.e., $X_i < u$) and $d-k$ variables where we are looking at $X_j > u$. Due to the perfect positive dependence ($X_i = X \sim U[0,1]$), the joint probability is
$$\mathbb{P}\left[X < \min_{i \in I} u \text{ and } X > \max_{j \in J} u\right] = \max\left(0, \min_{i \in I} u - \max_{j \in J} u\right)=0.$$
The integral $I_d(\alpha)$ just involves the integral of the product of the marginal probabilities. Therefore, we get:
$$I_d(\alpha) = -\frac{k! (d-k)!}{(d+1)!} .$$
Multiplying by the normalization factor, we obtain
$$\varphi_d^\alpha(M_d) = -\frac{2(d+1)}{(d-1)}\cdot \frac{k! (d-k)!}{(d+1)!}=-\frac{2}{(d-1)\binom{d}{k}}.$$
\item $k$ is odd ($0 < k < d$). When $k$ is odd, the joint probability is $$\mathbb{P}\left[\bigcap_{i=1}^d(\alpha_iX>\alpha_iu)\right] = \max\left(0, \min_{j \in J} u - \max_{i \in I} u\right)=0,$$
where $|I| = k$ and $|J| = d-k$. The formula for $\varphi_d^\alpha(M_d)$ is the same as the one derived for the even case of $k$. The difference in the directional dependence for odd $k$ arises from the structure of the joint probability involving $\min_{j \in J} u - \max_{i \in I} u$, which effectively reverses the roles of the minimum and maximum compared to the even $k$ case where we had $\min_{i \in I} u - \max_{j \in J} u$. However, the integral is equal to zero, leading to the same formula for $I_d(\alpha)$ and consequently for $\varphi_d^\alpha(M_d)$.
\end{enumerate}
In summary, where $k$ is the number of $-1$ in $\alpha$:
$$
\varphi_d^\alpha(M_d) = \begin{cases}
1 & \text{if } k = 0 \text{ or } k = d, \\
-\displaystyle\frac{2}{(d-1)\binom{d}{k}} & \text{otherwise}.
\end{cases}
$$
\end{example}

\begin{example}[Directional $\varphi$-coefficients for the Farlie-Gumbel-Morgenstern (FGM) family of $d$-copulas] \label{ex:FGM}Let us consider the FGM family of $d$-copulas given by
$$C_{\lambda}^{FGM}({\bf u})=\prod_{i=1}^du_i\left(1+\lambda\prod_{i=1}^d(1-u_i)\right),\quad\lambda\in[-1,1]$$
for all $\uu\in[0,1]^d$ (see \cite{Durante2016book,Ne06}). On the other hand, we have
$$(C_\lambda^{FGM})^\alpha(\uu)=\prod_{i=1}^du_i\left(1+(-1)^{|J|}\lambda\prod_{i=1}^d(1-u_i)\right),$$
    for all $\uu\in[0,1]^d$, where $J$ is the set of coordinates of $\alpha$ equal to $1$. Therefore,
    \begin{align*}
\varphi_d^\alpha(C_\lambda^{FGM})&=\frac{2(d+1)}{d-1}\int_0^1(-1)^{|J|}\lambda\prod_{i=1}^d\mathbb{P}[\alpha_iU_i>\alpha_iu](1-\mathbb{P}[\alpha_iU_i>\alpha_iu])du\\
        &=\frac{2(d+1)}{d-1}(-1)^{|J|}\lambda\int_0^1u^{|I|}(1-u)^{|I|}u^{|J|}(1-u)^{|J|}du\\
        &=\frac{2(d+1)}{d-1}(-1)^{|\mathcal{J}|}\lambda\int_0^1u^d(1-u)^ddu\\
        &=\frac{2(d+1)}{d-1}(-1)^{|\mathcal{J}|}\lambda\frac{(d!)^2}{(2d+1)!}\\
        &=\frac{2\lambda(-1)^{|J|}(d+1)(d!)^2}{(d-1)(2d+1)!}.
    \end{align*}
    So, for an even number of coordinates equal to 1 of $\alpha$, the coefficient takes the value 
     $$\varphi_d^\alpha(C_\lambda^{FGM})=\frac{2\lambda(d+1)(d!)^2}{(d-1)(2d+1)!},$$
     and for an odd number 
      $$\varphi_d^\alpha(C_\lambda^{FGM})=-\frac{2\lambda(d+1)(d!)^2}{(d-1)(2d+1)!}.$$

%\begin{figure}[h!]
%    \centering
%    \includegraphics[width=0.5\linewidth]
%{Figure_1.png}
%        \caption{Diagonal sections of the FGM %copulas associated with the different %directional random vectors for $d=3$ and %$\lambda=0.9$.}
%        \label{fig:placeholder}
%    \end{figure}
\end{example}

\begin{example}[Clayton copulas]
Let $C_\theta^{CL}$ be a Clayton $d$-copula, defined as
\begin{equation}\label{eq:Clayton}
   C_\theta^{CL}(\uu)=\left(\sum_{i=1}^du_i^{-\theta}-d+1\right)^{-1/\theta}, 
\end{equation}
for all $\uu\in[0,1]^d$ and $\theta>0$. Due to the symmetry of the copula it is not important which of the variables is transformed by changing the sign; only the number of signs changed matters. The general expression for the directional footrule coefficients for $C_\theta^{CL}$, using \eqref{th:expr}, is
$$\varphi_d^\alpha(C_\theta^{CL})=\frac{2(d+1)}{d-1}\sum_{k=0}^{|J|}(-1)^k\binom{|J|}{k}\int_0^1((|I|+k)u^{-\theta}+1-|I|-k)^{-1/\theta}\,du-\frac{2}{(d-1)\binom{d}{|J|}}.$$
\begin{figure}[h!]
        \centering
        \includegraphics[width=0.6\linewidth]{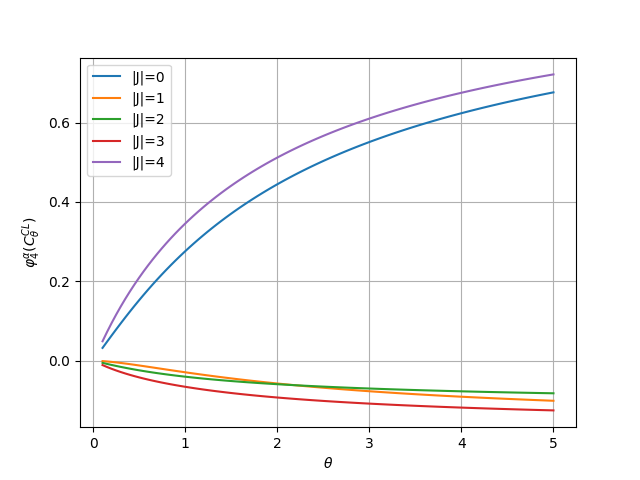}
        \caption{Values of $\varphi_4^\alpha(C_\theta^{CL})$ for Clayton $4$-copulas as functions of the parameter $\theta$, for all possible values of $|J|$. The curves for $|J|=0$ and $|J|=4$ (corresponding to $\alpha=-\mathbf{1}$ and $\alpha=\mathbf{1}$, respectively) exhibit positive dependence that increases with $\theta$, while intermediate values of $|J|$ display negative directional dependence.}
        \label{fig:placeholder}
    \end{figure}
\end{example}

\begin{example}[Cuadras-Augé $d$-copulas]
Let $C_\theta^{CA}$ be a multivariate Cuadras-Augé $d$-copula whose expression is given by
\begin{equation}\label{ej:CA}
    C_\theta^{CA}(\uu)=(\Pi(\uu))^{1-\theta}(M(\uu))^\theta,
\end{equation}
for all $\uu\in[0,1]^d$ and $\theta\in[0,1]$  (see \cite{Durante2016book,Ne06}). These $d$-copulas satisfy that $C_0^{CA}=\Pi$ and $C_1^{CA}=M$. Moreover, every lower $j$-dimensional marginal of this $d$-copula is a Cuadras-Augé $j$-copula. Since this copula is symmetric, it does not matter which components we keep equal to $1$, but their number. The general expression for these coefficients for $C_\theta^{CA}$ is, using \eqref{th:expr},
\begin{equation}\label{ex: expr}
    \varphi_d^\alpha\left(C_\theta^{CA}\right)=\frac{2(d+1)}{d-1}\sum_{k=0}^{|J|}(-1)^k\binom{|J|}{k}\frac{\theta(|I|+k-1)}{(|I|+k+1)^2-\theta[(|I|+k)^2-1]}.
\end{equation}

By using (\ref{ex: expr}), we obtain the value of $\varphi_4^\alpha(C_\theta^{CA})$ for all the possible components of $\alpha_i=1$
\begin{itemize}
    \item [i)] $|J|=0$ $$\varphi_4^-(C_\theta^{CA})=\frac{2\theta}{5-3\theta}.$$
    \item [ii)] $|J|=1$ $$\varphi_4^\alpha(C_\theta^{CA})=\frac{10}{3}\left(\frac{\theta}{8-4\theta}-\frac{3\theta}{25-15\theta}\right).$$
    \item [iii)] $|J|=2$ $$\varphi_4^\alpha(C_\theta^{CA})=\frac{10}{3}\left(\frac{\theta}{9-3\theta}-\frac{\theta}{4-2\theta}+\frac{3\theta}{25-15\theta}\right).$$
    \item [iv)] $|J|=3$ $$\varphi_4^\alpha(C_\theta^{CA})=\frac{10}{3}\left(-\frac{\theta}{3-\theta}+\frac{3\theta}{8-4\theta}-\frac{3\theta}{25-15\theta}\right).$$
    \item [v)] $|J|=4$ $$\varphi_4^\alpha(C_\theta^{CA})=\frac{10}{3} \left[\frac{2\theta}{3-\theta} - \frac{\theta}{2-\theta} + \frac{3\theta}{25-15\theta} \right].$$
\end{itemize}
In Figure \ref{fig:ex CA}
\begin{figure}[h!]
    \centering
    \includegraphics[width=0.6\linewidth]{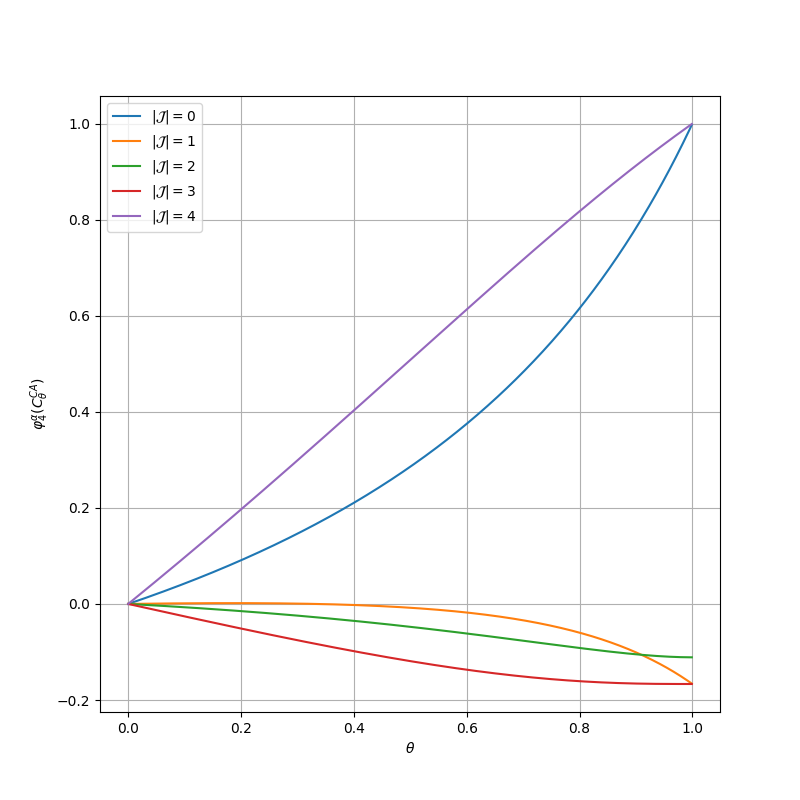}
    \caption{Values of $\varphi_4^{\alpha}$ for Cuadras-Augé $4$-copulas as a function of $\theta$, for all possible values of $\alpha$.}
    \label{fig:ex CA}
\end{figure}
we observe that, as anticipated, increasing values of the parameter $\theta$ drive the copula away from independence and toward maximal dependence in the cases $\alpha = \mathbf{1}$ and $\alpha = -\mathbf{1}$. When $\alpha$ has an odd number of its coordinates equal to $1$, we observe that the random variables show negative dependence.
\end{example}

The following result can be trivially extracted from Definition \ref{def:PDorder} and  the expression of the directional footrule-coefficients given by \eqref{eq:main2}.
\begin{corollary}
        Let $\XX$ and $\YY$ be two $d$-dimensional random vectors of continuous random variables uniformly distributed on $[0,1]$, whose respective joint distribution functions are the $d$-copulas $C_{\XX}$ and $C_{\YY}$. Let $\alpha=(\alpha_1,\ldots,\alpha_d)\in\mathbb{R}^d$ such that $|\alpha_i|=1$ for all $i=1,\ldots,d$. If $\XX\leq_{{\rm PD}(\alpha)}\YY$, then $\varphi_d^\alpha(C_{\XX})\leq\varphi_d^\alpha(C_{\YY})$.\\
\end{corollary}

As an illustrative case, let $C_{\lambda_1}^{FGM}$ and $C_{\lambda_2}^{FGM}$ denote two members of the FGM family of $d$-copulas, as described in Example \ref{ex:FGM}. It is well established that these copulas are ordered with respect to the ${\rm PD}(\alpha)$ ordering, namely, $
C_{\lambda_1}^{FGM} \leq_{{\rm PD}(\alpha)} C_{\lambda_2}^{FGM}$
if, and only if, their associated parameters satisfy $\lambda_1 \leq \lambda_2$ (see \cite{Edeamo2024}). As a direct implication of this ordering property, the correlation measure $\varphi_d^\alpha$ is monotone with respect to the parameter $\lambda$, which yields
$
\varphi_d^\alpha\bigl(C_{\lambda_1}^{FGM}\bigr) \leq \varphi_d^\alpha\bigl(C_{\lambda_2}^{FGM}\bigr)
$
if, and only if, $\lambda_1 \leq \lambda_2$.

\section{Nonparametric estimators}\label{sec:esti}

Our goal now is to estimate the directional coefficients $\varphi_d^\alpha$. As shown in (\ref{eq:esp}), directional footrule-coefficients can be regarded in terms of expectations. Several rank-based estimators have been proposed in the literature (see, e.g., \cite{Edeamo25}, \cite{Dec25}, \cite{Pe2016}, \cite{Pe2023}); following this line of work, we introduce a nonparametric rank-based estimator for the directional version of Spearman's footrule. As shown in \cite{Pe2016}, the use of nonparametric rank-based estimators offers several advantages, as alternative methods—such as estimators based on the empirical copula—may produce values that fall outside the parametric space. Consider a $d$-dimensional random sample of size $n$, $\{X_{1j},\ldots,X_{dj}\}_{j=1}^n$, of a $d$-dimensional random vector $(X_1,\ldots,X_d)$ with associated $d$-copula $C$. Let $R_{ij}$ be the rank of $X_{ij}$ in $\{X_{i1},\ldots,X_{in}\}$. Then we define, from (\ref{eq:esp}), and using the pseudo-observations $\widehat{U}_{ij}=R_{ij}/(n+1)$
$$\widetilde{\varphi}_{n,d}^\alpha=\frac{2(d+1)}{(d-1)(n+1)}\left\{\frac{1}{n}\sum_{j=1}^n\left(\min_{i\in J} R_{ij}-\max_{i\in I}R_{ij}\right)_+-\frac{(n+1)|I|!|J|!}{(d+1)!}\right\},$$
where $A_+=\max\{0,A\}$. Since $\mathbb{E}\left[\left(\min_{i\in J}U_i-\max_{i\in I}U_i\right)_+\right]$ can be estimated by 
$$\frac{1}{n}\sum_{j=1}^n\left(\min_{i\in J}\widehat{U}_{ij}-\max_{i\in I}\widehat{U}_{ij}\right)_+=\frac{1}{n(n+1)}\sum_{j=1}^n\left(\min_{i\in J}R_{ij}-\max_{i\in I}R_{ij}\right)_+.$$

Note that when $\alpha=\mathbf{1}$ and $\alpha=-\mathbf{1}$, it holds that 
$$\widetilde{\varphi}_{n,d}^+=\frac{2(d+1)}{(d-1)(n+1)}\left(\frac{1}{n}\sum_{j=1}^n\min_{1\leq i\leq d}R_{ij}-\frac{n+1}{d+1}\right)\quad\text{and}\quad \widetilde{\varphi}_{n,d}^-=\frac{2(d+1)}{(d-1)(n+1)}\left(\frac{1}{n}\sum_{j=1}^n\min_{1\leq i\leq d}\overline{R}_{ij}-\frac{n+1}{d+1}\right),$$
since when $J=\emptyset$, the term $\min_{i\in J}R_{ij}-\max_{i\in I}R_{ij}$ turns to be $n+1-\max_{1\leq i\leq d}R_{ij}$, or equivalently, $\min_{1\leq i\leq d}\overline{R}_{ij}$. The expressions obtained for the extreme cases of the directional footrule estimators are consistent and coherent with the previously derived expressions presented in \cite{Dec25}. 

It can be shown that the proposed estimators admit a decomposition as a linear combination of the lower-dimensional $\widetilde{\varphi}$-estimators, following a structure analogous to that established in Theorem \ref{th:exp}:

\begin{theorem}{\label{th:est}}
    Let $\{X_{1j},\ldots,X_{dj}\}_{j=1}^n$ be a $d$-dimensional random sample of size $n$ of a $d$-dimensional random vector $(X_1,\ldots,X_d)$ with associated $d$-copula $C$. Let $\alpha\in\mathbb{R}^d$ be such that $\alpha_i\in\{-1,1\}$ for all $i=1,\ldots,d$. Let $I\subseteq\{1,2,\ldots,d\}$ such that $\alpha_i=-1$ if $i\in I$ and $\alpha_i=1$ if $i\in J=\{1,2,\ldots,d\}\setminus I$. Then 
    \begin{equation}
        \widetilde{\varphi}_{n,d}^\alpha=\frac{2(d+1)}{d-1}\sum_{K\subseteq J}(-1)^{|K|}\frac{|I|+|K|-1}{2(|I|+|K|+1)}\widetilde{\varphi}^-_{n,I\cup K}.
    \end{equation}
\end{theorem}

\begin{proof}
Since $\min_{i\in L}\overline{R}_{ij}=n+1-\max_{i\in L}R_{ij}$ for every index subset $L$, we can rewrite the expression of the estimators as follows
    \begin{align*}
        \widetilde{\varphi}_{n,d}^\alpha&=\frac{2(d+1)}{(d-1)(n+1)}\left\{\frac{1}{n}\sum_{j=1}^n\left(\min_{i\in J} R_{ij}-\max_{i\in I}R_{ij}\right)_+-\frac{(n+1)|I|!|J|!}{(d+1)!}\right\}\\
        &=\frac{2(d+1)}{(d-1)(n+1)}\left\{\frac{1}{n}\sum_{j=1}^n\left(\min_{i\in I} \overline{R}_{ij}-\max_{i\in J}\overline{R}_{ij}\right)_+-\frac{(n+1)|I|!|J|!}{(d+1)!}\right\}.
    \end{align*}
    Note that $(\min_{i \in I}\overline{R}_{ij} - \max_{i \in J} \overline{R}_{ij})_+ = 0$ when $\min_{i \in I} \overline{R}_{ij} < \max_{i \in J} \overline{R}_{ij}$. Therefore, in the following we can restrict our attention to the case where $\min_{i \in I} \overline{R}_{ij} \geq \max_{i \in J} \overline{R}_{ij}$. We consider now that $\min_{i\in I}\overline{R}_{ij}\geq \max_{i\in J}\overline{R}_{ij}$; otherwise, the argument equals zero. Taking this into account, it holds that
    \begin{align*}
        \widetilde{\varphi}_{n,d}^\alpha&=\frac{2(d+1)}{(d-1)(n+1)}\left\{\frac{1}{n}\sum_{j=1}^n\left(\min_{i\in I} \overline{R}_{ij}-\max_{i\in J}\overline{R}_{ij}\right)_+-\frac{(n+1)|I|!|J|!}{(d+1)!}\right\}\\
        &=\frac{2(d+1)}{(d-1)(n+1)}\left\{\frac{1}{n}\sum_{j=1}^n\min_{i\in I} \overline{R}_{ij}-\frac{1}{n}\sum_{j=1}^n\max_{i\in J}\overline{R}_{ij}-\frac{(n+1)|I|!|J|!}{(d+1)!}\right\}\\
        &=\frac{2(d+1)}{(d-1)(n+1)}\left\{\frac{1}{n}\sum_{j=1}^n\min_{i\in I} \overline{R}_{ij}-\sum_{k=1}^{|J|}(-1)^{k+1}\sum_{\substack{K\subseteq J\\ |K|=k}}\frac{1}{n}\sum_{j=1}^n\min_{i\in K}\overline{R}_{ij}-\frac{(n+1)|I|!|J|!}{(d+1)!}\right\},
    \end{align*}
    where the following inclusion-exclusion expression has been used
    $$\max(x_1,...,x_d)=\sum_{k=1}^d(-1)^{k+1}\sum_{\substack{K\subseteq{1,...,d}\\|K|=k}}\min_{i\in K}x_i.$$
   Since  $\min_{i\in I}\overline{R}_{ij}\geq \max_{i\in J}\overline{R}_{ij}$, then $\min_{i\in I\cup K}\overline{R}_{ij}=\min(\min_{i\in I}\overline{R}_{ij},\min_{i\in K}\overline{R}_{ij})=\min_{i\in K}\overline{R}_{ij}$. Thus
   \begin{align*}
       \widetilde{\varphi}_{n,d}^\alpha&=\frac{2(d+1)}{(d-1)(n+1)}\left\{\frac{1}{n}\sum_{j=1}^n\min_{i\in I} \overline{R}_{ij}-\sum_{k=1}^{|J|}(-1)^{k+1}\sum_{\substack{K\subseteq J\\ |K|=k}}\frac{1}{n}\sum_{j=1}^n\min_{i\in K}\overline{R}_{ij}-\frac{(n+1)|I|!|J|!}{(d+1)!}\right\}\\
       &=\frac{2(d+1)}{(d-1)(n+1)}\left\{\sum_{k=0}^{|J|}(-1)^k\sum_{\substack{K\subseteq J\\ |K|=k}}\frac{1}{n}\sum_{j=1}^n\min_{i\in I\cup K}\overline{R}_{ij}-\frac{(n+1)|I|!|J|!}{(d+1)!}\right\}\\
       &=\frac{2(d+1)}{(d-1)(n+1)}\left\{\sum_{K\subseteq J}\frac{(-1)^{|K|}}{n}\sum_{j=1}^n\min_{i\in I\cup K}\overline{R}_{ij}-\frac{(n+1)|I|!|J|!}{(d+1)!}\right\}\\
       &=\frac{2(d+1)}{(d-1)(n+1)}\left\{\sum_{K\subseteq J}\left((-1)^{|K|}\frac{(|I|+|K|-1)(n+1)}{2(|I|+|K|+1)}\widetilde{\varphi}^-_{n,I\cup K}+\frac{n+1}{|I|+|K|+1}\right)-\frac{(n+1)|I|!|J|!}{(d+1)!}\right\}\\
       &=\frac{2(d+1)}{(d-1)}\left\{\sum_{K\subseteq J}\left((-1)^{|K|}\frac{(|I|+|K|-1)}{2(|I|+|K|+1)}\widetilde{\varphi}^-_{n,I\cup K}+\frac{1}{|I|+|K|+1}\right)-\frac{|I|!|J|!}{(d+1)!}\right\}.
   \end{align*}
Computing the terms in the sum that do not involve the estimators, we obtain 
   $$\sum_{K\subseteq J}(-1)^{|K|}\frac{1}{|I|+|K|+1}=\frac{|I|!|J|!}{(d+1)!}.$$
Hence, this term cancels the last term appearing in the expression, yielding the desired linear combination of the estimators.
\end{proof}

\begin{figure}[h!]
        \centering
        \includegraphics[width=0.8\linewidth]{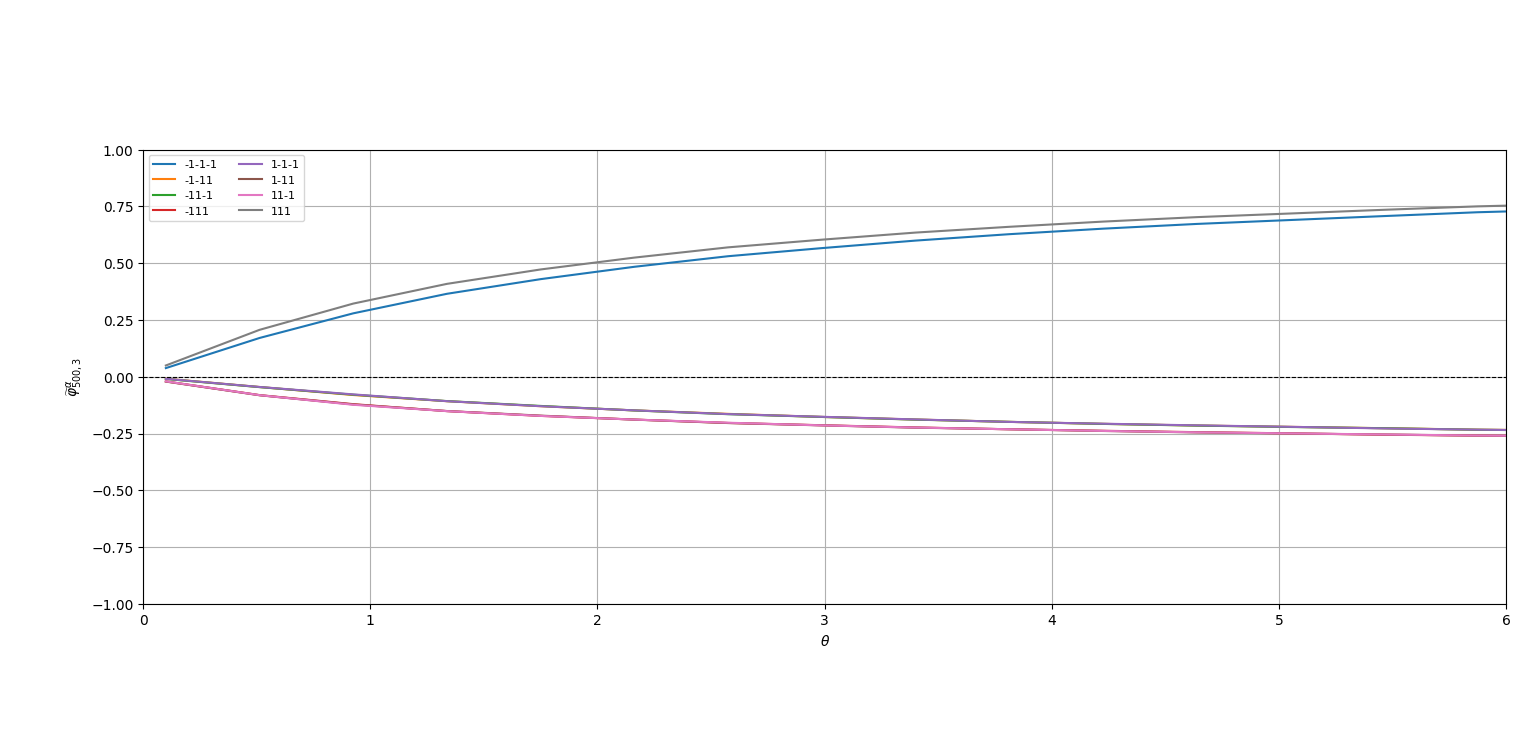}
        \caption{Values of the different estimators of the directional footrule-coefficients as functions of the parameter $\theta$ for $3$-dimensional Clayton family copulas.}
        \label{fig:clayton3D}
    \end{figure}

We now illustrate the good performance of these estimators through a practical example in the $3$-dimensional setting. Table \ref{tab:1} presents the $3$-dimensional directional footrule-coefficients for all possible directions $\alpha\in\{-1,1\}^3$, computed exactly via numerical integration and subsequently approximated using our proposed estimators, based on $1000$ independent samples, each one of size $500$. 
    \begin{table}[h!]
        \centering
        \begin{tabular}{c c c }
           $\alpha$  & $\varphi_d^\alpha$ & $\widetilde{\varphi}_{n,d}^\alpha$ \\ \hline
         $(-1,-1,-1)$ & $0.68975$ & $0.68924$\\
         $(-1,1,1)$ & $-0.24935$ & $-0.24848$\\
        $(1,-1,1)$ & $-0.24935$ &  $-0.24889$ \\
        $(1,1,-1)$ & $-0.24935$ & $-0.24874$\\
        $(1,-1,-1)$ & $-0.22020$ & $-0.22004$\\
        $(-1,-1,1)$ & $-0.22020$ & $-0.22030$\\
        $(-1,1,-1)$ & $-0.22020$ & $-0.22045$\\
        $(1,1,1)$ & $0.71890$ & $0.71768$ 
        \end{tabular}
        \caption{Results obtained in Monte Carlo simulations in order to calculate the value of $\widetilde{\varphi}^\alpha_{n,d}$ in every direction $\alpha\in\{-1,1\}^3$ based on a total of 1000 samples of size 500 generated for the Clayton 3-copula for $\theta=5$.}
        \label{tab:1}
    \end{table}
To further illustrate the behavior of the proposed estimators, Figure \ref{fig:clayton3D} presents their values as functions of parameter $\theta$ for three-dimensional Clayton family copulas. The figure highlights how the estimators capture the variation in directional dependence across different values of $\theta$, providing a visual confirmation of their finite-sample performance and their alignment with the theoretical directional footrule-coefficients.

\subsection{Empirical processes and some properties of these estimators}

We now establish the asymptotic properties of the proposed directional footrule estimators. Building on the theory of empirical copula processes, we derive conditions under which these estimators are consistent and asymptotically normal. This theoretical foundation enables statistical inference in practical applications and validates the use of these estimators for directional dependence analysis.

\subsubsection{Empirical process}

In \cite{Gar2013}, an empirical process was introduced in order to estimate a generalization of a $3$-copula which can be easily extended to the $n$-dimensional case. As the authors did, we commence considering $\upbeta$ an index set, such that $\upbeta\subseteq\{1,\ldots,d\}$. We define $x_\upbeta=(x_{\upbeta_1},\ldots,x_{\upbeta_k})$ an arbitrary value of $\mathbf{X}_\upbeta=(X_i)_{i\in\upbeta}$. Let $|\upbeta|$ denote the cardinal of $\upbeta$ and $\alpha\in\mathbb{R}^d$ a direction. Consider the function 
    \begin{equation}\label{H}
     H_{\upbeta,\alpha}(x_\upbeta)=\mathbb{P}\left[\alpha_iX_i\leq\alpha_ix_i, \text{ for all } i\in\upbeta\right],
    \end{equation}
i.e.,
$$
H_{\upbeta,\alpha}(x_\upbeta) = \mathbb{P}\left[ \bigcap_{\substack{i \in \upbeta \\ \alpha_i=1}} \{X_i \leq x_i\}, \bigcap_{\substack{i \in \upbeta \\ \alpha_i=-1}} \{X_i \geq x_i\} \right].
$$

Let $u_\upbeta$ be the analogous of $x_\upbeta$ and $F_i$ the marginal cumulative distribution function of $X_i$. Then 
\begin{equation}\label{C}
C_{\upbeta,\alpha}(u_\upbeta) = H_{\upbeta,\alpha}((v_i)_{i \in \upbeta}),
\end{equation}
where
$$
v_i =
\begin{cases}
F_i^{-1}(u_i), & \text{if } \alpha_i = 1, \\
F_i^{-1}(1-u_i), & \text{if } \alpha_i = -1,
\end{cases}$$
is a generalization of the $d$-copula when $\alpha_i=1$ and $\upbeta=\{1,\ldots,d\}$. From this expression, we can introduce, as in \cite{Gar2013}, the empirical process to estimate it
\begin{equation}\label{estC}
C_{\upbeta,\alpha,n}(u_\upbeta)=\frac{1}{n+1}\sum_{j=1}^n\prod_{i\in\upbeta}{\bf 1}_{\left\{\alpha_i\frac{R_{ij}}{n+1}\leq\alpha_iu_i\right\}}.    
\end{equation}
When $\upbeta=\{1,\ldots,d\}$ and $\alpha={\bf 1}$, Expression (\ref{estC}) is the empirical estimator of the $d$-copula.

To determine the weak convergence ---denoted by ``$\xrightarrow[]{w}$''--- of the empirical process 
$$\sqrt{n}\{C_{\upbeta,\alpha,n}(u_\upbeta)-C_{\upbeta,\alpha}(u_\upbeta)\}, \hspace{0.5cm} u_{\upbeta}\in [0,1]^{|\upbeta|},$$
it suffices to proceed as in \cite{Gar2013}, where the following two conditions are considered:

\begin{condition}\label{con1}
    For each $i\in\upbeta$, the $i$-th first-order partial derivative $\dot{C}_{i,\upbeta,\alpha}$ exists and is continuous on the set $\{u_\upbeta\in[0,1]^{|\upbeta|}:0<u_i<1\}$.
\end{condition}

\begin{condition}\label{con2}
$\mathbb{B}_{C_{\upbeta,\alpha}}(u_\upbeta)$ is a $C_{\upbeta,\alpha}$-tight centered Gaussian process on $[0,1]^{|\upbeta|}$. $u_{\upbeta}^{(i)}$ is a vector such that its $j$-th component, for $j \in \upbeta$, is:
$$(u_{\upbeta}^{(i)})_j=\begin{cases}
1 & \text{if } j \neq i \text{ when } \alpha_j = 1, \\
0 & \text{if } j \neq i \text{ when } \alpha_j = -1, \\
u_i & \text{if } j = i.
\end{cases}$$
The covariance function is $\text{Cov}(\mathbb{B}_{C_{\upbeta,\alpha}}(u_\upbeta),\mathbb{B}_{C_{\upbeta,\alpha}}(v_\upbeta))=C_{\upbeta,\alpha}(w_\upbeta)-C_{\upbeta,\alpha}(u_\upbeta)C_{\upbeta,\alpha}(v_\upbeta)$, where the $j$-th component $(w_\upbeta)_j$ is:
$$(w_\upbeta)_j=\begin{cases}
\min(u_j,v_j) & \text{if } \alpha_j = 1, \\
\max(u_j,v_j) & \text{if } \alpha_j = -1,
\end{cases}$$
for all $j \in \upbeta$.
\end{condition}

From Conditions \ref{con1} and \ref{con2}, the following result is obtained, which is valid for any dimension $d$ and, therefore, guarantees the weak convergence of our process.

\begin{theorem}\label{th:proc}
    Let $H_{\upbeta,\alpha}$ be a $|\upbeta|$-dimensional distribution function, given by (\ref{H}) with continuous marginal distributions $F_i$, $i\in\upbeta$ and with $C_{\upbeta,\alpha}$ given by (\ref{C}), where $\upbeta\subseteq\{1,...,d\}$ and $\alpha\in\mathbb{R}^d$, $\alpha_i\in\{-1,1\}$. Under Condition \ref{con1} on the function $C_{\upbeta,\alpha}$, when $n\longrightarrow\infty$
    $$\sqrt{n}\left\{C_{\upbeta,\alpha,n}(u_\upbeta)-C_{\upbeta,\alpha}(u_\upbeta)\right\}\xrightarrow[]{w}\mathbb{G}_{C_{\upbeta,\alpha}}(u_\upbeta).$$
    Weak convergence takes place in $l^\infty([0,1]^{|\upbeta|})$ and $\mathbb{G}_{C_{\upbeta,\alpha}}(u_\upbeta)=\mathbb{B}_{C_{\upbeta,\alpha}}(u_\upbeta)-\sum_{i\in\upbeta}\dot{C}_{i,\upbeta,\alpha}(u_\upbeta)\mathbb{B}_{C_{\upbeta,\alpha}}(u_\upbeta^{(i)})$, where $\mathbb{B}_{C_{\upbeta,\alpha}}$ is a $C_{\upbeta,\alpha}$-tight centered Gaussian process on $[0,1]^{|\upbeta|}$ and $\mathbb{G}_{C_{\upbeta,\alpha}}$ follows Condition \ref{con2}.
\end{theorem}

\subsubsection{Properties of the estimators}\label{sect:prop}
Our main objective in this section is to prove that the estimators we have defined are unbiased. To do so, first, we must emphasize that our estimators verify the following:
\begin{equation}\label{eq:concl}
   \widetilde{\varphi}_{n,d}^\alpha=\frac{2(d+1)(n+1)}{(d-1)n}\int_0^1C_{\upbeta,-\alpha,n}(u,\ldots,u)\,{\rm d}u-\frac{2}{(d-1)\binom{d}{|J|}},
\end{equation}
since $$\int_0^1C_{\upbeta,-\alpha,n}(u,\ldots,u)\,{\rm d}u=\int_0^1\frac{1}{n+1}\sum_{j=1}^n\prod_{i=1}^d\mathbf{1}_{\left\{\alpha_i\frac{R_{ij}}{n+1}\geq \alpha_iu\right\}}\,{\rm d}u=\frac{1}{(n+1)^2}\sum_{j=1}^n\left(\min_{i\in J}R_{ij}-\max_{i\in I}R_{ij}\right)_+$$

In order to prove that our estimators are unbiased we will use the following theorem from \cite{Gar2013}, which is an adaptation of \cite[Theorem 6]{Fer2004} and \cite{Gan87}.

\begin{theorem}\label{th:proc2}
    Under the assumptions of Theorem \ref{th:proc}, let  $\{a_n\}_{n\geq 1}$ and $\{b_n\}_{n\geq 1}$ be two sequences of real numbers satisfying $\sqrt{n}(a_n-a_0)=O(n^{-1/2})$ and $\sqrt{n}(b_n-b_0)=O(n^{-1/2})$, respectively, where $a_0$ and $b_0$ are constant values. Let $$T_n(f)=a_n\int_{I^{|\upbeta|}}f(u_\upbeta){\rm d}u_\upbeta+b_n,$$
    for $n\geq 0$, where $f$ is a $|\upbeta|$-integrable function. Then, when $n\longrightarrow\infty$,
    $$\sqrt{n}\{T_n(C_{\upbeta,\alpha,n})-T_0(C_{\upbeta,\alpha})\}\xrightarrow{w}Z_{C_{\upbeta,\alpha}}\sim N(0,\sigma_{C_{\upbeta,\alpha}}^2),$$
    with $$\sigma_{C_{\upbeta,\alpha}}^2=a_0^2\int_{I^{|\upbeta|}}\int_{I^{|\upbeta|}}\mathbb{E}[\mathbb{G}_{C_{\upbeta,\alpha}}(u_\upbeta)\mathbb{G}_{C_{\upbeta,\alpha}}(v_\upbeta)]{\rm d}u_\upbeta {\rm d}v_\upbeta$$ and $$Z_{C_{\upbeta,\alpha}}:=a_0\int_{I^{|\upbeta|}}\mathbb{G}_{C_{\upbeta,\alpha}}(u_\upbeta){\rm d}u_\upbeta.$$
\end{theorem}

It is clear that this theorem remains valid if we take the $b_n$ to be constant in \cite[Proof of Theorem 5.2]{Gar2013} . As a consequence, we have the following result.

\begin{corollary}\label{col:unb}
    Under assumptions of Theorem \ref{th:proc}, when $n\longrightarrow\infty$
    $$\sqrt{n}\{\widetilde{\varphi}_{n,d}^\alpha-\varphi_d^\alpha\}\xrightarrow{w}Z_{C_{\upbeta,-\alpha}}\sim N(0,\sigma_{C_{\upbeta,-\alpha}}^2),$$
    where $\upbeta=\{1,2,\ldots,d\}$ and $\alpha\in\mathbb{R}^d$ is a vector such that $\alpha_i\in\{-1,1\}$ for all $i=1,\ldots,d$.
\end{corollary}

\begin{proof}
    We consider the real number sequences in Theorem \ref{th:proc2} as $$a_n=\frac{2(d+1)(n+1)}{(d-1)n}\quad\text{and}\quad a_0=\frac{2(d+1)}{d-1}.$$
    We are going to check that the condition $\sqrt{n}(a_n-a_0)=O(n^{-1/2})$ holds. It is straightforward since
    \begin{align*}
        \sqrt{n}(a_n-a_0)=&\sqrt{n}\left(\frac{2(d+1)(n+1)}{(d-1)n}-\frac{2(d+1)}{d-1}\right)\\
        &=\sqrt{n}\left(\frac{2(d+1)}{(d-1)n}\right)\\
        &=O(n^{-1/2}).
    \end{align*}

    So, using (\ref{eq:concl}), the desired result is obtained.
\end{proof}

By Corollary \ref{col:unb}, the estimators $\widetilde{\varphi}_{n,d}^\alpha$ are asymptotically unbiased and Var$(\widetilde{\varphi}_{n,d}^\alpha)\longrightarrow 0$ whenever $n\longrightarrow\infty$. Therefore, $\widetilde{\varphi}_{n,d}^\alpha$ is also an asymptotically consistent estimator since convergence in probability, $\widetilde{\varphi}_{n,d}^\alpha\xrightarrow{p}\varphi_d^\alpha$, is guaranteed  as a consequence of Chebyshev's inequality.

With the results obtained above, the conjecture posed in \cite{Dec25}, involving the extreme footrule-estimators $\widetilde{\varphi}_{n,d}^+$ and $\widetilde{\varphi}_{n,d}^-$, is now resolved. Specifically, it is established that the proposed estimators are consistent and asymptotically normally distributed, provided that certain regularity conditions are satisfied. This confirmation not only validates the theoretical properties anticipated in the earlier work, but also provides a rigorous foundation for the application of these estimators in practical settings, including multivariate dependence modeling and directional coefficient estimation.

\subsection{Simulations}

In order to evaluate the empirical behavior of the proposed estimators, Monte Carlo simulations for $d$-copulas of the parametric Clayton and Cuadras-Augé families have been carried out.  The expressions of these copulas are given by (\ref{eq:Clayton}) and (\ref{ej:CA}), respectively. Note that both the Clayton and Cuadras-Aug\'e copulas satisfy Conditions \ref{con1} and \ref{con2} due to their smoothness properties on the interior of the unit hypercube. Condition \ref{con1} holds since both families possess continuous partial derivatives, while Condition \ref{con2} is guaranteed by standard empirical copula process theory \cite{Fer2004,Gan87}, ensuring the asymptotic validity of the proposed estimators.

Two dimensions $d=\{4,5\}$ have been considered and four sample sizes $n=\{20,50,100,500\}$. A Monte Carlo study was conducted in order to evaluate the behavior of the proposed estimators under different sampling scenarios. For each combination of direction and parameter value, $1000$ independent samples of size $n$ were generated. Based on each sample, the directional estimators $\widetilde{\varphi}_{n,d}^{\alpha}$ were computed, and their empirical averages were subsequently obtained. These averages are reported in the corresponding tables, together with reference values of the associated directional coefficients $\varphi_d^{\alpha}$, obtained through numerical integration.

The simulation results are organized according to the dimension, direction, and copula parameter. In particular, Table~\ref{tab:3} summarizes the estimated values of $\varphi_4^{(-1,1,1,-1)}$ for the $4$-dimensional Clayton copula, considering sample sizes $n \in \{20,50,100,500\}$ and parameter values $\theta \in \{0.4, 0.6, 1, 2, 5\}$. Table~\ref{tab:4} presents the estimation results for the coefficient $\varphi_5^{(1,-1,1,-1,1)}$ associated with the five-dimensional Clayton copula, again for the same range of sample sizes and parameter values. In all cases, the numerical approximations of the true coefficients serve as benchmarks for assessing the accuracy of the proposed estimators.

\begin{table}[h!]
\begin{center}
\begin{tabular}{ c  c  c  c  c c c  }
\hline
\multicolumn{3}{ c }{Sample size} & $n=20$  & $n=50$ & $n=100$ & $n=500$ \\ \hline
$\theta$ & $\varphi_4^{\alpha}$ & $\widetilde{\varphi}_{n,4}^{\alpha}$ & mean & mean & mean & mean \\ \hline
0.4 & -0.02042 & $\widetilde{\varphi}_{n,4}^{(-1,1,1,-1)}$ & -0.02660 & -0.02147 & -0.02198  & -0.02036 \\
0.6 & -0.02829 &  $\widetilde{\varphi}_{n,4}^{(-1,1,1,-1)}$ & -0.03366  & -0.02891 & -0.02960  & -0.02816  \\
1 & -0.04053  &  $\widetilde{\varphi}_{n,4}^{(-1,1,1,-1)}$ & -0.04402 & -0.04090 & -0.04141 & -0.04029 \\
2 & -0.05946  & $\widetilde{\varphi}_{n,4}^{(-1,1,1,-1)}$ & -0.06104  & -0.05920 & -0.06000 & -0.05915 \\
5 & -0.08234  &  $\widetilde{\varphi}_{n,4}^{(-1,1,1,-1)}$ & -0.08166 &  -0.08144 & -0.08218 & -0.08213 \\ \hline
\end{tabular}
\caption{Results obtained in Monte Carlo simulations in order to calculate the value of $\widetilde{\varphi}_{n,4}^{(-1,1,1,-1)}$ based on a total of $1000$ samples of size $n$ generated for the Clayton $4$-copula for different values of the parameter $\theta$.}
\label{tab:3}
\end{center}
\end{table}

\begin{table}[h!]
\begin{center}
\begin{tabular}{ c  c  c  c  c c c  }
\hline
\multicolumn{3}{ c }{Sample size} & $n=20$  & $n=50$ & $n=100$ & $n=500$ \\ \hline
$\theta$ & $\varphi_5^{\alpha}$ & $\widetilde{\varphi}_{n,5}^{\alpha}$ & mean & mean & mean & mean \\ \hline
0.4 & -0.01105 & $\widetilde{\varphi}_{n,5}^{(1,-1,1,-1,1)}$ & -0.01144 & -0.01246 & -0.01064  & -0.01125 \\
0.6 & -0.01483 &  $\widetilde{\varphi}_{n,5}^{(1,-1,1,-1,1)}$ & -0.01607  & -0.01475 & -0.01491  & -0.01475  \\
1 & -0.02042  &  $\widetilde{\varphi}_{n,5}^{(1,-1,1,-1,1)}$ & -0.02135 & -0.02085 &  -0.02135 &- 0.02055 \\
2 & -0.02864  & $\widetilde{\varphi}_{n,5}^{(1,-1,1,-1,1)}$ & -0.02897  & -0.02892 & -0.02863 & -0.02865 \\
5 & -0.03821  &  $\widetilde{\varphi}_{n,5}^{(1,-1,1,-1,1)}$ & -0.03756 &  -0.03789 & -0.03749 & -0.03817 \\ \hline
\end{tabular}
\caption{Results obtained in Monte Carlo simulations in order to calculate the value of $\widetilde{\varphi}_{n,5}^{(1,-1,1,-1,1)}$ based on a total of $1000$ samples of size $n$ generated for the Clayton $5$-copula for different values of the parameter $\theta$.}
\label{tab:4}
\end{center}
\end{table}

A careful inspection of the numerical results reveals that the proposed estimators closely track the corresponding theoretical directional footrule-coefficients, with their accuracy improving as the sample size increases. Across all configurations examined, the estimated coefficients take negative values for the directions under consideration. Nevertheless, for the Clayton $d$-copula, extreme directions lead to positive coefficient values, which increase toward unity as the dependence parameter $\theta$ grows, in agreement with the expected dependence behavior of the model.

Table \ref{tab:6} displays the satisfactory performance of the estimators in all admissible directions in $\mathbb{R}^4$ for $4$-copulas from the Cuadras--Augé family, considering in particular two copulas associated with specific parameter values, $\theta=\{0.4,0.8\}$. The results are based on Monte Carlo simulations conducted on 1000 samples of different sizes, $n=\{20,50,100,500\}$, with the table reporting the outcomes for the case $n=500$.

\begin{table}[h!]
    \centering
    % Primera subtabla
    \begin{subtable}[h]{0.45\textwidth}
        \centering
        \begin{tabular}{p{2cm} c c}
            Parameter & $\alpha$ & $\widetilde{\varphi}_{n,4}^\alpha$ \\ \hline
            \multirow{16}{*}{$\theta=0.4$} 
            & (-1,-1,-1,-1) & 0.38348 \\
            & (1,-1,-1,-1) & -0.06336 \\
            & (-1,1,-1,-1) & -0.063197 \\
            & (1,1,-1,-1) & -0.04292 \\
            & (-1,-1,1,-1) & -0.06376 \\
            & (1,-1,1,-1) & -0.04352 \\
            & (-1,1,1,-1) & -0.04279 \\
            & (1,1,1,-1) & -0.06392 \\
            & (-1,-1,-1,1) & -0.06372 \\
            & (1,-1,-1,1) & -0.04319 \\
            & (-1,1,-1,1) & -0.04314 \\
            & (1,1,-1,1) & -0.06395 \\
            & (-1,-1,1,1) & -0.04230 \\
            & (1,-1,1,1) & -0.06362 \\
            & (-1,1,1,1) & -0.06458 \\
            & (1,1,1,1) & 0.38450 \\
        \end{tabular}
    \end{subtable}
    % Segunda subtabla
    \begin{subtable}[h]{0.45\textwidth}
        \centering
        \begin{tabular}{p{2cm} c c}
            Parameter & $\alpha$ & $\widetilde{\varphi}_{n,4}^\alpha$ \\ \hline
            \multirow{16}{*}{$\theta=0.8$} 
            & (-1,-1,-1,-1) & 0.77980 \\
            & (1,-1,-1,-1) & -0.12969 \\
            & (-1,1,-1,-1) & -0.12967 \\
            & (1,1,-1,-1) & -0.08708 \\
            & (-1,-1,1,-1) & -0.12984 \\
            & (1,-1,1,-1) & -0.08672 \\
            & (-1,1,1,-1) & -0.08692 \\
            & (1,1,1,-1) & -0.12987 \\
            & (-1,-1,-1,1) & -0.12956 \\
            & (1,-1,-1,1) & -0.08715 \\
            & (-1,1,-1,1) & -0.08731 \\
            & (1,1,-1,1) & -0.12933 \\
            & (-1,-1,1,1) & -0.08718 \\
            & (1,-1,1,1) & -0.12965 \\
            & (-1,1,1,1) & -0.12932 \\
            & (1,1,1,1) & 0.77950 \\
        \end{tabular}
    \end{subtable}
    \caption{Results obtained in Monte Carlo simulations in order to calculate the value of $\widetilde{\varphi}^\alpha_{n,4}$ in every direction $\alpha\in\{-1,1\}^4$ based on a total of 1000 samples of size 500 generated for the Cuadras-Augé 4-copula for $\theta=0.4$ and $\theta=0.8$.}
    \label{tab:6}
\end{table}

\begin{figure}[h!]
    \centering
    \begin{subfigure}[h]{0.45\linewidth}
    \includegraphics[width=\linewidth]{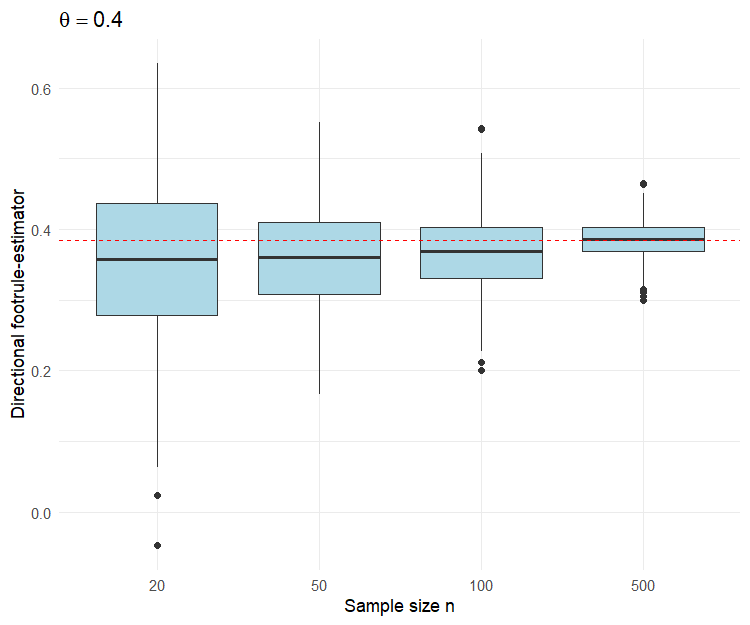}
\end{subfigure}
    \begin{subfigure}[h]{0.45\linewidth}
    \includegraphics[width=\linewidth]{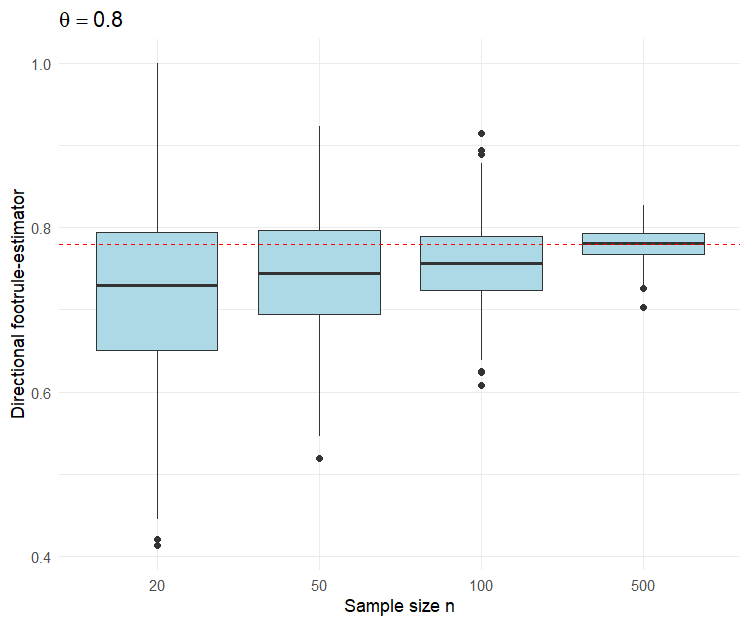}
\end{subfigure}
\caption{Empirical distribution of $\widetilde{\varphi}_4^{(1,1,1,1)}$ for the Cuadras-Augé $4$-copula with parameter values $\theta=\{0.4,0.8\}$ and sample sizes $n=\{20,50,100,500\}$. Horizontal red line shows the value of $\varphi_4^{(1,1,1,1)}$.}
\label{fig: diagbig1}
\end{figure}

\begin{figure}[h!]
    \centering
    \begin{subfigure}[h]{0.45\linewidth}
    \includegraphics[width=\linewidth]{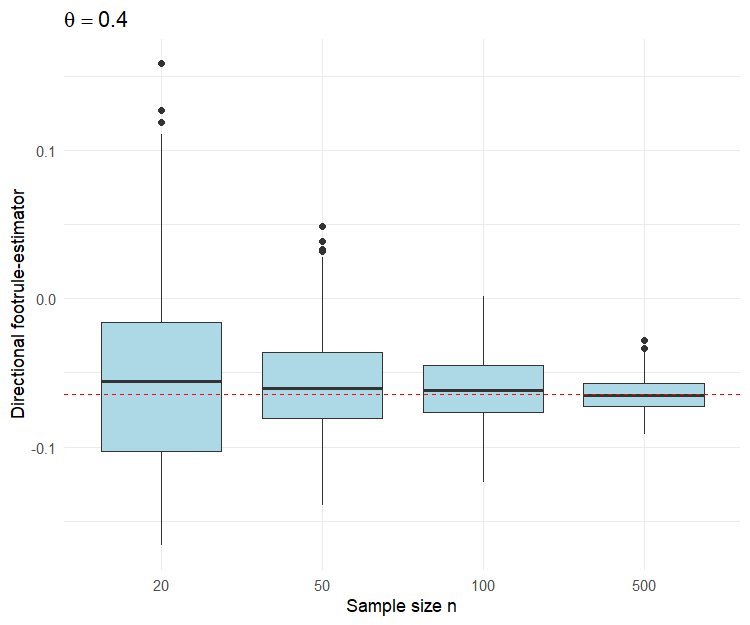}
\end{subfigure}
    \begin{subfigure}[h]{0.45\linewidth}
    \includegraphics[width=\linewidth]{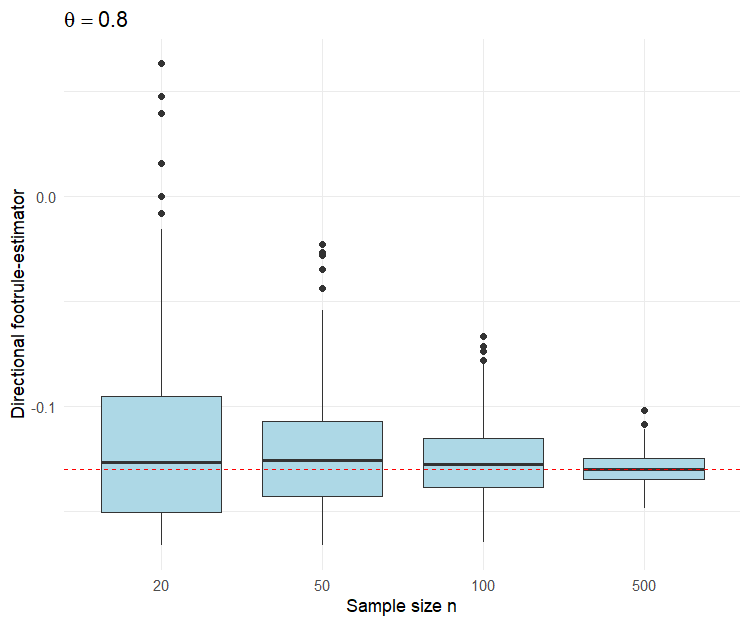}
\end{subfigure}
\caption{Empirical distribution of $\widetilde{\varphi}_4^{(1,-1,1,1)}$ for the Cuadras-Augé $4$-copula with parameter values $\theta=\{0.4,0.8\}$ and sample sizes $n=\{20,50,100,500\}$. Horizontal red line shows the value of $\varphi_4^{(1,-1,1,1)}$.}
\label{fig: diagrambig2}
\end{figure}

The analysis of this additional example once again highlights the strong performance of the proposed estimators when applied to the Cuadras--Augé family of copulas. The numerical results confirm their ability to accurately capture the underlying dependence structure across different directional configurations. Furthermore, as illustrated in Figures \ref{fig: diagbig1} and \ref{fig: diagrambig2}, the estimators exhibit clear asymptotic consistency, with convergence toward the corresponding theoretical values as the sample size increases.

These properties make the proposed directional footrule-estimators particularly suitable for practical applications involving multivariate dependence modeling, such as risk management, reliability analysis, and financial or environmental data. In such contexts, Cuadras--Augé copulas provide a flexible framework for modeling asymmetric dependence structures, and the reliability of the estimators ensures their applicability within a broad class of mathematical models where directional dependence plays a relevant role.

\begin{figure}[h!]
    \centering
    \includegraphics[width=0.6\linewidth]{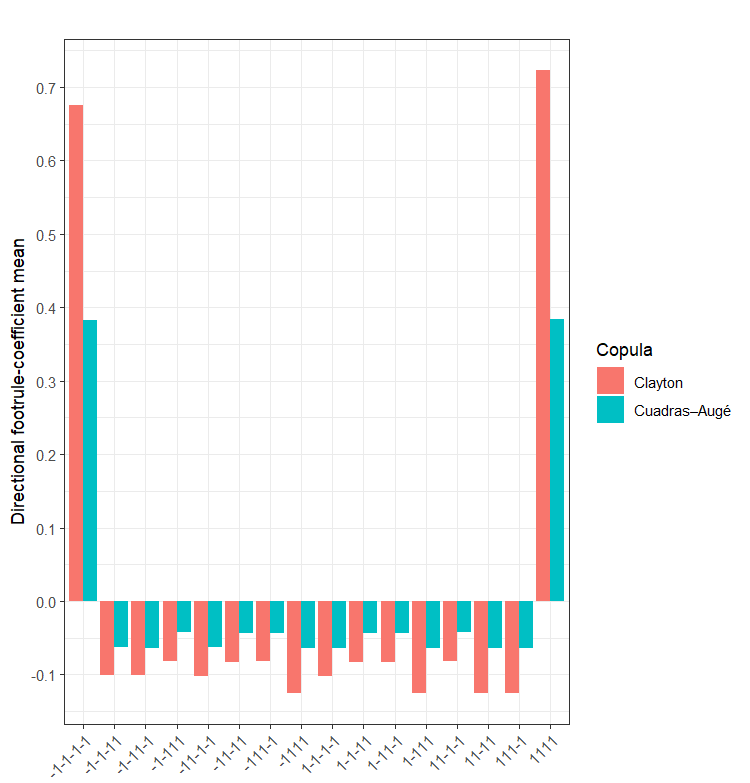}
    \caption{Comparative study of the estimated directional footrule-coefficients for two $4$-dimensional copulas corresponding to the parametric families of Clayton and Cuadras-Augé. For the Clayton Copula we have considered a parameter value $\theta=5$, and for the Cuadras-Auge, $\theta=0.4$. }
    \label{fig:vs}
\end{figure}

Finally, to conclude this analysis, Figure~\ref{fig:vs} presents a comparative plot of the estimated directional footrule-coefficients for both the Clayton and Cuadras--Augé copulas. In the case of the four-dimensional Clayton copula, the directions $\mathbf{1}$ and $-\mathbf{1}$ exhibit positive coefficients. Simulations have made the symmetry of the copula evident, as reflected by the equality of the coefficients depending on the number of coordinates equal to $-1$

In contrast, for the four-dimensional Cuadras--Augé copula, the direction $(1,1,1,1)$ clearly dominates, with an estimated coefficient close to the dependence parameter $\theta$ (approximately 0.4 in our example). Mixed directions containing both $+1$ and $-1$ components yield coefficients that are essentially zero, indicating that the mass of the distribution is concentrated along the comonotonic diagonal. These results illustrate the singular-plus-continuous nature of the Cuadras--Augé copula and highlight the ability of the proposed estimators to capture subtle differences in dependence structure across copula families

\section{Conclusions}\label{sec:conc}

In this paper, we have introduced a new family of directional multivariate footrule-coefficients within the copula framework, aimed at capturing direction-specific dependence structures in multivariate data. By extending the classical Spearman's footrule to a directional setting, the proposed coefficients provide a refined description of dependence patterns that cannot be detected by traditional measures.

We have established the main theoretical properties of the proposed coefficients, showing their consistency with existing multivariate versions in the extreme directional cases and highlighting their structural coherence. In particular, the representation of the directional coefficients as linear combinations of lower-dimensional directional footrule-coefficients offers new insight into the contribution of partial dependence among subsets of components to the overall dependence structure.

To facilitate practical implementation, we have proposed nonparametric rank-based estimators for the directional coefficients and studied their finite-sample behavior through Monte Carlo simulations. The numerical results demonstrate that the proposed estimators provide accurate approximations of the theoretical coefficients and exhibit improved performance as the sample size increases. Explicit expressions and illustrative examples for well-known copula families, including the Farlie--Gumbel--Morgenstern and Cuadras--Aug\'e copulas, further emphasize the interpretability and flexibility of the proposed approach.

Overall, the directional footrule-coefficients introduced in this work enrich the toolbox of rank-based dependence measures and offer a versatile framework for analyzing complex multivariate dependence structures. Possible directions for future research include the study of their asymptotic distributions, the development of statistical tests for directional dependence, and the extension of the proposed methodology to broader classes of dependence models.

\section{Acknowledgments}

The first author is partially supported by the CDTIME of the University of Almería (Spain).


\begin{thebibliography}{99}

\bibitem{Ahn2015} Ahn, J.Y. (2015). Negative dependence concept in copulas and the marginal free herd behavior index. {\it J. Comput. Appl. Math.} {\bf 288}, 304--322.

\bibitem{Edeamo} Amo, E. de, García-Fernández, D., Quesada-Molina, J.J., Úbeda-Flores, M. (2025). Modeling directional monotonicity with copulas. {\it Iranian J. Fuzzy Syst.} {\bf 22}, 135--146.

\bibitem{Edeamo25} Amo, E. de, García-Fernández, D., Úbeda-Flores, M. (2026). Directional $\varphi$-coefficients. To appear.

\bibitem{Edeamo2024} Amo, E. de, Rodríguez-Griñolo, M.R., Úbeda-Flores, M. (2024). Directional dependence orders of random vectors. {\it Mathematics }{\bf 12}, article 419.

\bibitem{Chen} Chen, C., Xu, W., Zhang, W. Zhu, H., Dai, D. (2023). Asymptotic properties of Spearman's footrule and Gini's gamma in bivariate normal model. {\it J. Frankl. Inst.} {\bf 360}(13), 9812--9843.

\bibitem{Dec20} Decancq, K. (2020). Measuring cumulative deprivation and affluence based on the diagonal dependence diagram. {\it Metron} {\bf 78}, 103--117.

\bibitem{Dec25} Decancq, K., Pérez, A., Prieto-Alaiz, M. (2025). Multivariate Dependence Based on Diagonal Sections: Spearman’s Footrule and Related Measures. In: {\it Steland, A., Rafajłowicz, E., Parolya, N. (eds) Stochastic Models, Statistics and Their Applications}. SMSA 2024. Springer Proceedings in Mathematics \& Statistics, vol 499. Springer, Cham.

\bibitem{Dol06} Dolati, A., Úbeda-Flores, M. (2006) On mesaures of multivariate concordance. {\it J. Probab. Stat. Sci.}{\bf 4}(2), 147--163.

\bibitem{Durante2016book} Durante, F., Sempi, C. (2016). {\it Principles of Copula Theory}. Chapman $\&$ Hall/CRC, Boca Raton.

 \bibitem{Fer2004} Fermanian, J.D., Radulovic, D., Wegkamp, M. (2004). Weak convergence of empirical copula processes. {\it Bernoulli} {\bf 10}, 847--860.

\bibitem{Fisher97} Fisher, N.I. (1997). Copulas. In: {\it Encyclopedia of statistical sciences}, Vol. 1 (S. Kotz, C.B. Read, D.L. Banks, Eds.), Wiley, New York, pp. 159--163.

\bibitem{Fuc14} Fuchs, S. (2014). Multivariate copulas: Transformations, symmetry, order and measures of concordance. {\it Kybernetika} {\bf 50}(5), 725--743.

\bibitem{Gan87} Gaenssler, P., Stute, W. (1987). {\it Seminar on Empirical Processes}. Springer Basel, Basel.

\bibitem{Gar2013} García, J.E., González-López, V.A., Nelsen, R.B. (2013). A new index to measure positive dependence in trivariate distributions. {\it J. Multivariate Anal.} {\bf 115}, 481--495. 

\bibitem{Genest} Genest, C., Nešlehová, J., Ghorbal, N.B. (2010). Spearman's footrule and Gini's gamma: a review with complements. {\it J. Nonparam. Stat.} {\bf 22}, 937--954.

\bibitem{Graham} Graham, R.L., Knuth, D.E., Patashnik, O. (1994). \emph{Concrete Mathematics (2nd ed.)}. Addison-Wesley, reading.

\bibitem{Jodgeo82} Jogdeo, K. (1982). Concepts of dependence. In: {\it Encyclopedia of Statistical Sciences}, Vol. 1 (S. Kotz, N.L. Johnson, Eds.), Wiley, New York, pp. 324--334.

\bibitem{Joe1990} Joe, H. (1990). Multivariate concordance. {\it J. Multivariate Anal.} {\bf 35}, 12--30.

\bibitem{Joe1997} Joe, H. (2014). {\it Dependence modeling with copulas}. Chapman $\&$ Hall, New York.

\bibitem{Mu06} M\"uller, A., Scarsini, M. (2006). Archimedean copulae and positive dependence. {\it J. Multivariate Anal.} {\bf 93}, 434--445.

%\bibitem{Navarro2021} Navarro, J., Pellerey, F., Sordo, M.A. 
%(2021). Weak dependence notions and their mutual %relationships. {\em Mathematics} {\bf 9}, article 81.

\bibitem{Ne96} Nelsen, R.B. (1996). Nonparametric measures of multivariate association. In: {\it Distributions with Fixed Marginals and Related Topics}, Vol. 28 (L. Rüschendorf, B. Schweizer, M.D. Taylor, Eds.), Institute of Mathematical Statistics, Hayward, 223--232.

\bibitem{Ne06} Nelsen, R.B. (2006). {\it An Introduction to Copulas (2nd ed.)}. Springer, New York.

\bibitem{Ne2011} Nelsen, R.B., Úbeda-Flores, M. (2011). Directional dependence in multivariate distributions. {\it Ann. Inst. Stat. Math} {\bf 64}, 677--685. 

\bibitem{Pe2016} Pérez, A., Prieto-Alaiz, M. (2016). A note on nonparametric estimation of copula-based multivariate extensions of Spearman's rho. {\em Stat. Probab. Lett.} {\bf 112}, 41--50.

\bibitem{Pe2023}Pérez, A. Prieto-Alaiz, M., Chamizo, F., Liebscher, E., Úbeda-Flores, M. (2023). Nonparametric estimation of the multivariate Spearman's footrule: A further discussion. {\it Fuzzy Sets and Systems} {\bf 467}, 108489.

\bibitem{Que2012} Quesada-Molina, J.J., Úbeda-Flores, M. (2012). Directional dependence of random vectors. {\it Inf. Sci.} {\bf 215}, 67--74.

\bibitem{Que2024} Quesada-Molina, J.J., Úbeda-Flores, M. (2024). Monotonic random variables according to a direction. {\em Axioms} {\bf 13}, 275.

\bibitem{Sk59} Sklar, A. (1959). Fonctions de r\'epartition $\grave{\rm a}$ $n$ dimensions et leurs marges. {\it Publ. Inst. Statist. Univ. Paris} {\bf 8}, 229--231.

\bibitem{Sp06} Spearman, C. (1906). 'Footrule' for measuring correlation. {\it British Journal of Psychology} {\bf 2}, 89--108.

\bibitem{Ub05} Úbeda-Flores, M. (2005) Multivariate versions of Blomqvist's beta and Spearman's footrule. {\it Ann. Inst. Statist. Math.} {\bf 57}, 781--788.

\bibitem{Ub2017} \'Ubeda-Flores, M., Fern\'andez-S\'anchez, J. (2017). Sklar's theorem: The cornerstone of the Theory of Copulas. In: {\it Copulas and Dependence Models with Applications} (M. \'Ubeda-Flores, E. de Amo Artero, F. Durante, J. Fern\'andez-S\'anchez, Eds.), Springer, Cham, pp. 241--258.

\bibitem{Wei2014} Wei, Z., Wang, T., Panichkitkosolkul, W. (2014). Dependence and association concepts through copulas. In: {\it Modeling Dependence in Econometrics - Advances in Intelligent Systems and Computing}, Vol. 251 (V.N. Huynh, V. Kreinovich, S. Sriboonchitta, Eds.), Springer, Cham, pp. 113--126.

\bibitem{Wolff1980} Wolff, E.F. (1980). $N$-dimensional measures of dependence. {\it Stochastica} {\bf 4(3}), 175--188.

\end{thebibliography}
\end{document}